\documentclass[12pt]{article}
\usepackage{amssymb}
\usepackage{amsfonts}
\usepackage{latexsym}
\usepackage{amsthm}

\theoremstyle{plain}
\theoremstyle{remark}

\begin{document}
\def\a{\alpha}  \def\cA{{\cal A}}     \def\bA{{\bf A}}  \def\mA{{\mathscr A}}
\def\b{\beta}   \def\cB{{\cal B}}     \def\bB{{\bf B}}  \def\mB{{\mathscr B}}
\def\g{\gamma}  \def\cC{{\cal C}}     \def\bC{{\bf C}}  \def\mC{{\mathscr C}}
\def\G{\Gamma}  \def\cD{{\cal D}}     \def\bD{{\bf D}}  \def\mD{{\mathscr D}}
\def\d{\delta}  \def\cE{{\cal E}}     \def\bE{{\bf E}}  \def\mE{{\mathscr E}}
\def\D{\Delta}  \def\cF{{\cal F}}     \def\bF{{\bf F}}  \def\mF{{\mathscr F}}
\def\c{\chi}    \def\cG{{\cal G}}     \def\bG{{\bf G}}  \def\mG{{\mathscr G}}
\def\z{\zeta}   \def\cH{{\cal H}}     \def\bH{{\bf H}}  \def\mH{{\mathscr H}}
\def\e{\eta}    \def\cI{{\cal I}}     \def\bI{{\bf I}}  \def\mI{{\mathscr I}}
\def\p{\psi}    \def\cJ{{\cal J}}     \def\bJ{{\bf J}}  \def\mJ{{\mathscr J}}
\def\vT{\Theta} \def\cK{{\cal K}}     \def\bK{{\bf K}}  \def\mK{{\mathscr K}}
\def\k{\kappa}  \def\cL{{\cal L}}     \def\bL{{\bf L}}  \def\mL{{\mathscr L}}
\def\l{\lambda} \def\cM{{\cal M}}     \def\bM{{\bf M}}  \def\mM{{\mathscr M}}
\def\L{\Lambda} \def\cN{{\cal N}}     \def\bN{{\bf N}}  \def\mN{{\mathscr N}}
\def\m{\mu}     \def\cO{{\cal O}}     \def\bO{{\bf O}}  \def\mO{{\mathscr O}}
\def\n{\nu}     \def\cP{{\cal P}}     \def\bP{{\bf P}}  \def\mP{{\mathscr P}}
\def\r{\rho}    \def\cQ{{\cal Q}}     \def\bQ{{\bf Q}}  \def\mQ{{\mathscr Q}}
\def\s{\sigma}  \def\cR{{\cal R}}     \def\bR{{\bf R}}  \def\mR{{\mathscr R}}
\def\S{\Sigma}  \def\cS{{\cal S}}     \def\bS{{\bf S}}  \def\mS{{\mathscr S}}
\def\t{\tau}    \def\cT{{\cal T}}     \def\bT{{\bf T}}  \def\mT{{\mathscr T}}
\def\f{\phi}    \def\cU{{\cal U}}     \def\bU{{\bf U}}  \def\mU{{\mathscr U}}
\def\F{\Phi}    \def\cV{{\cal V}}     \def\bV{{\bf V}}  \def\mV{{\mathscr V}}
\def\P{\Psi}    \def\cW{{\cal W}}     \def\bW{{\bf W}}  \def\mW{{\mathscr W}}
\def\o{\omega}  \def\cX{{\cal X}}     \def\bX{{\bf X}}  \def\mX{{\mathscr X}}
\def\x{\xi}     \def\cY{{\cal Y}}     \def\bY{{\bf Y}}  \def\mY{{\mathscr Y}}
\def\X{\Xi}     \def\cZ{{\cal Z}}     \def\bZ{{\bf Z}}  \def\mZ{{\mathscr Z}}
\def\O{\Omega}
\def\ve{\varepsilon}
\def\vt{\vartheta}
\def\vp{\varphi}

\def\Z{{\Bbb Z}}
\def\R{{\Bbb R}}
\def\C{{\Bbb C}}
\def\T{{\Bbb T}}
\def\N{{\Bbb N}}
\def\S{{\Bbb S}}

\def\ma{\left(\begin{array}{cc}}
\def\am{\end{array}\right)}
\def\iint{\int\!\!\!\int}
\def\lt{\biggl}
\def\rt{\biggr}
\let\ge\geqslant
\let\le\leqslant
\def\[{\begin{equation}}
\def\]{\end{equation}}
\def\wt{\widetilde}
\def\pa{\partial}
\def\sm{\setminus}
\def\es{\emptyset}
\def\no{\noindent}
\def\ol{\overline}
\def\iy{\infty}
\def\ev{\equiv}
\def\/{\over}
\def\ts{\times}
\def\os{\oplus}
\def\ss{\subset}
\def\h{\hat}
\def\na{\nabla}
\def\Re{\mathop{\rm Re}\nolimits}
\def\Im{\mathop{\rm Im}\nolimits}
\def\supp{\mathop{\rm supp}\nolimits}
\def\sign{\mathop{\rm sign}\nolimits}
\def\Ran{\mathop{\rm Ran}\nolimits}
\def\Ker{\mathop{\rm Ker}\nolimits}
\def\Tr{\mathop{\rm Tr}\nolimits}
\def\const{\mathop{\rm const}\nolimits}
\def\Wr{\mathop{\rm Wr}\nolimits}
\def\da{\downarrow}
\def\BBox{\hspace{1mm}\vrule height6pt width5.5pt depth0pt \hspace{6pt}}
%\def\liml{\lim\limits}
%\def\iintl{\iint\limits}
%\def\bl{\biggl}

%%%%%%%%%%%%%%%

\def\Twelve{
\font\Tenmsa=msam10 scaled 1200 \font\Sevenmsa=msam7 scaled 1200
\font\Fivemsa=msam5 scaled 1200
%\newfam\msafam
\textfont\msafam=\Tenmsa \scriptfont\msafam=\Sevenmsa
\scriptscriptfont\msafam=\Fivemsa

\font\Tenmsb=msbm10 scaled 1200 \font\Sevenmsb=msbm7 scaled 1200
\font\Fivemsb=msbm5 scaled 1200
%\newfam\msafam
\textfont\msbfam=\Tenmsb \scriptfont\msbfam=\Sevenmsb
\scriptscriptfont\msbfam=\Fivemsb

\font\Teneufm=eufm10 scaled 1200 \font\Seveneufm=eufm7 scaled 1200
\font\Fiveeufm=eufm5 scaled 1200
%\newfam\eufmfam
\textfont\eufmfam=\Teneufm \scriptfont\eufmfam=\Seveneufm
\scriptscriptfont\eufmfam=\Fiveeufm}

\def\Ten{
\textfont\msafam=\tenmsa \scriptfont\msafam=\sevenmsa
\scriptscriptfont\msafam=\fivemsa

\textfont\msbfam=\tenmsb \scriptfont\msbfam=\sevenmsb
\scriptscriptfont\msbfam=\fivemsb

\textfont\eufmfam=\teneufm \scriptfont\eufmfam=\seveneufm
\scriptscriptfont\eufmfam=\fiveeufm}

\title
{Schr\"odinger operator with a  junction of two 1-dimensional
periodic potentials}

\author{Evgeny Korotyaev
\begin{footnote}  {Institut f\"ur  Mathematik, Humboldt Universit\"at zu Berlin,
Rudower Chaussee 25, 12489, Berlin, Germany,\ \ \ e-mail:
evgeny@math.hu-berlin.de }
\end{footnote}
}

\maketitle

\begin{abstract}

\no The spectral properties of the Schr\"odinger operator $T_ty=
-y''+q_ty$ in $L^2(\R )$ are studied, with a potential
$q_t(x)=p_1(x), x<0, $ and  $q_t(x)=p(x+t), x>0, $ where $p_1, p$ are
periodic potentials and $t\in \R$ is a parameter of dislocation.
Under some conditions  there exist simultaneously gaps in the
continuous spectrum of $T_0$ and eigenvalues in these gaps. The
main goal of this  paper is to study the discrete spectrum and the
resonances of $T_t$. The following results are obtained: i) In any
gap of $T_t$ there exist $0,1$ or $2$ eigenvalues. Potentials with
0,1 or 2 eigenvalues in the gap are constructed. ii) The
dislocation, i.e. the case $p_1=p$ is studied. If $t\to 0$, then in
any gap in the spectrum there exist both eigenvalues ($ \le 2 $) and
resonances ($ \le 2 $) of $T_t$ which belong to a gap on the
second sheet and  their asymptotics as $t\to 0 $ are determined.
iii) The eigenvalues of the half-solid, i.e. $p_1={\rm constant}$,
are also studied. iv) We prove that for any even 1-periodic
potential $p$ and any sequences $\{d_n\}_1^{\iy }$, where $d_n=1$
or $d_n=0$  there exists a unique  even 1-periodic potential $p_1$
with the same gaps and $d_n$ eigenvalues of $T_0$ in the n-th gap
for each $n\ge 1.$

\end{abstract}

\section {Introduction}
\setcounter{equation}{0}

%   \hskip 1cm

We consider the Schr\"odinger  operator $T_ty=-y''+q_t(x)y$ acting on $L^2(\R)$,
where the real potential $q_t$  is given by
\[
q_t(x)=\cases {p_1(x)\ \ \ &if\ \ \ \ $x<0,\ \ \ \ p_1\in L^1(\t\T) $\cr
          p(x+t)\ \ \ &if\ \ \ \ $x>0, \ \ \ \  p\in L^1(\T)$\cr}, \ \ \
          t\in [0,1],
\]
where $\t>0$ is the period of $p_1$ and $\T=\R/\Z$. We call such a potential $q_t$ biperiodic. If $p_1=p$ and $t=0,$ then we obtain  the well known periodic case,
that is  the Hill operator $H =-{d^2\/dx^2}+p$ in $L^2(\R )$.
It is well known (see [T]) that the spectrum of $H $ is absolutely
continuous and consists of intervals $  \sigma _n=[\a^+_{n-1}, \a^-_n ],$ where $\a^+_{n-1}< \a^-_n\le \a^+_{n},\ n\ge 1$. These intervals are separated by gaps $\g_n(H)=(\a^-_n, \a^+_n ), n\ge 1$. We set $\g _0(H)=(-\iy ,\a_0^+)$ and $\g(H)=\cup \g_n(H)$.
If a gap degenerates, i.e. $\g_n=\es, n\ge 1$, then the corresponding segments $\s_n,\s_{n+1}$ merge. Define the Hill operator $H(t) = -{d^2\/dx^2}+p(x+t)$
in $L^2(\R )$. It is clear that $\s (H)=\s (H(t))$ for any $t\in [0,1]$. Introduce the Hill operator $H_m =-{d ^2\/dx^2}+p_m(x)$ in $L^2(\R ), m=1, 2$ where here and below $p_2\ev p$.
In the biperiodic case (see Theorem 2.1) we prove that
the spectrum of each  $T_t, t\in \R$ has the following form
\[
\s (T_t)=\s _{ac}(T_t)\cup \s _{d}(T_t),
\ \ \ \  \s _{ac}(T_t)=\s (H_1)\cup \s (H_2) ,
\ \ \ \ \s _d(T_t)\ss \g (H_1)\cap \g (H_2).
\]
Note that  (1.2) implies $\s_{sc}(T_t)=\es$ and there are no
embedded eigenvalues. The absolutely continuous spectrum
$\s_{ac}(T_t)$ consists of intervals $\s_n(T_t), n\ge 1$. These
intervals are separated by the gaps $\g_n(T_t), n\ge 1$. In
general, there exist eigenvalues in these gaps and we have (see
Theorem 2.1)
\[
\# (T_t,\g_n(T_t))\le 2, \ \ \ \  {\rm for \ \ each }\ \ n\ge 0,
\]
where $\# (T_t,\o )$ is the number of eigenvalues of $T_t$ on an interval $\o $.
The basic goal of this paper is to study eigenvalues in these gaps.

In this paper we study in more detail two  cases of biperiodic
potentials. Firstly, we consider the case of dislocation,  that is
the operator $T_t^{di}= -{d^2\/dx^2}+p_{(t)}(x)$ in $L^2(\R ),$ with the potential $p_{(t)}(x)=\c_-(x)p(x)+\c_+(x)p(x+t), t\in \R$, here and below $\c_{\pm}(x)=1, \pm x\ge 0$ and
$\c_{\pm}(x)=0, \pm x<0$ and the potential $p$  as above. Changing $t$ we get different potentials $p_{(t)}(\cdot )$. It can be shown that the eigenvalues of $T_t^{di}$ are periodic in $t$ (see (1.4)). Actually, they are 2-periodic, not 1-periodic. 
If $t=0$, then $p_{(0)}(x)=p(x)$ is periodic and $T_0^{di}=H$ and
eigenvalues are absent. But if $t\neq 0$ what are we able to say
about the eigenvalues in the gaps $\g _n, n\ge 0,$ of the operator
$T_t^{di}$?

Secondly, we consider the half-solid, i. e. the operator $T_t^sy=-y''+q_t^s(x)y$ in $L^2(\R )$,  where $q_t^s(x)=s\c_-(x)+\c_+(x)p(x+t)$ and $t,s\in\R$.
Changing constants $t, s$ we get different potentials $q_t^s(x)$. The eigenvalues of $T_t^s$ depend on $t, s$.
If the first gap of $H$ is open and $s>\a_1^-$, then there exists a gap in the spectrum of $T_t^s$. What are we able to say on the eigenvalues in these gaps $\g _n(T_t^s), n\ge 0$?

We describe the main results which are proved in the present  paper:

\no  {\it i) relations (1.2-3),

\no ii) the description of the possible coexistence of eigenvalues and resonances for different classes of potentials,

\no iii) we prove that for any fixed numbers $m (m=0,1,2)$ and
$N\ge 1$ there exists a biperiodic potential such that $\# (T_0,\g_n(T_0))=m$ for all $1\le n\le N$,

\no iv) the discrete spectrum of $T_t, T_t^{di}, T_t^s$ is studied.
For example, if $t$ is small then there exist eigenvalues of $T_t^d$ in the gaps $\g _n, n\ge 0$, and their asymptotics are determined as $t\to 0$.

\no v) The analytic continuation of $(T_0-\l)^{-1}, \Im \l>0$, into the second sheet of the energy surface is obtained and the resonances are studied. For the dislocation, $T_t^{di}=T_t$, we have: $\l_1$ is an eigenvalue of $T_t$,
iff the number $\l_1$ considered as a point on the second sheet is a resonance of the operator $\wt T_t$ given by $\wt T_ty =-y''+\wt q_t(x)y$, where $\wt q_t(x)=\c_-(x)p_2(x+t)+\c_+(x)p_1(x)$.}

We briefly describe the proof. Using  results of self-adjoint
extensions [Kr] we obtain (1.2-3). An eigenvalue of $T_t$ is a
zero of the corresponding Wronskian. Using techniques from
inverse spectral theory [L], [Tr], [K2] we study how such a zero
depends on $t.$ Then we prove (1.3) with $\# (T,\g_n(T))=m, m=0,
1,$ for all $n\ge 1.$ In order to obtain the case m=2 we use
additional results from  [GT] (see also [KK1]) on the inverse
problem for the Hill operator. In order  to determine the
asymptotics of the eigenvalues as $t\to 0$
we study the corresponding Wronskian. Here we use the implicit function Theorem, various properties of the fundamental solutions and the quasimomentum.

 Tamm [Ta] was the first, who considered the Schr\"odinger operator with biperiodic potentials (for the Kronig-Penny model) and proved the existence of eigenvalues (the famous surface states). The spectral properties of the Schr\"odinger operator with biperiodic
potentials were studied in various papers (see e.g. [A] [A1],
[DS]).   The spectral problem with biperiodic potentials arises in non-linear equations. For example, Bikbaev and Sharipov [BS] considered the KdV equation with biperiodic  initial data. 
In this paper [BS] it is assumed that
eigenvalues are absent. Anoshchenko [A] studied the inverse
problem for the Schr\"odinger operator with a biperiodic potential
plus a decreasing one. But the author does not study 
the discrete spectrum of the Schr\"odinger operator
 with a biperiodic potential. The author [K3] obtained
the following result concerning the motion of eigenvalues
for the case of the dislocation:
in each gap $\g_n\neq \es, n\ge 1$ there exist two unique "states"
(an eigenvalue and a resonance) $\l_n^{\pm}(t)$ of the dislocation operator, such that $\l_n^{\pm}(0)=\a_n^{\pm}$ and
the point $\l_n^{\pm}(t)$ runs clockwise around the gap $\g_n$
 changing the energy sheet whenever it hits $\a_n^{\pm}$, making
$n/2$ complete revolutions in unit time.
On the first sheet $\l_n^{\pm}(t)$ is an eigenvalue
and on the second sheet $\l_n^{\pm}(t)$ is a resonance.
Moreover, the following identities are fulfilled:
\[
\l_{2n}^{\pm}(t+1)=\l_{2n}^{\pm}(t),\ \ and  \ \ \l_{2n+1}^{\pm}(t+1)=-\l_{2n+1}^{\pm}(t)
+\a_{2n+1}^-+\a_{2n+1}^+, \ \ 
any \ n\ge 0, t\in \R.
\] 
We think that the results of the present paper are needed to study the KdV equation, the inverse problem and the spectral properties of $T_t+V(x)$, where $V(x)\to 0, $ as $x\to \pm \iy $ or $V=\ve x$ ( Stark effect), and junctions of two d-dimensional periodic potentials, $d>1$.

\section {Main results}
\setcounter{equation}{0}

We recall some properties of the Hill operator $H=-{d^2\/dx^2}+p$
acting in $L^2(\R )$, where the potential $p\in L^1(\T)$
is real. The spectrum of $H $ is absolutely continuous  and consists of intervals $\s_n=[\a^+_{n-1},\a^-_n],$  where $\a^+_{n-1}<\a^-_n\le\a^+_{n}, n\ge 1$.
We set $\a^+_0=0.$ These intervals are separated by gaps $\g _n(H)=(\a^-_n, \a^+_n ), n\ge 1$. Let $\vp(x,z), \vt(x,z)$ be the solutions of the equation
\[
-y''+py=\l y, \ \ \ \l\in \C ,
\]
satisfying $\vp_x(0,\l)=\vt(0,\l)=1,$ and
$\vp(0,\l)=\vt_x(0,\l)=0$. We define the Lyapunov function
$\D(\l)={1\/2}(\vp_x(1,\l)+\vt (1,\l))$ and the function
$a(\l)={1\/2}(\vp_x(1,\l)-\vt (1,\l)).$ Note that
$\D(\a_{n}^{\pm})=(-1)^n,\ n\ge 1$. The sequence
$\a_0^+<\a_1^-\le\a_1^+\ <\dots$ is the spectrum of Eq. (2.1) with
2-periodic boundary conditions, that is $y(x+2)=y(x), x\in\R$. Let
$\P (x,\a_n^{\pm})$ be the corresponding real normalized
eigenfunctions, i.e. $\int _0^1\P (x, \a_n^{\pm})^2dx=1$. 
If $\a_n^-= \a_n^+$, then $\a_n^-$ is a double
eigenvalue. The lowest eigenvalue $\a_0^+$ is simple,
$\D(\a_0^+)=1,$ and the corresponding eigenfunction has period 1.
The eigenfunctions corresponding to $\a_n^{\pm}$ have period 1
when $n$ is even and they are anti-periodic, $\P
(x+1,\a_n^{\pm})=-\P (x,\a_n^{\pm}),\ x\in\R $, when $n$ is odd.
Let $\mu _n(p), n\ge 1, $ be the Dirichlet spectrum of (2.1) with
the boundary condition $y(0)=y(1)=0$. Let $\n_n(p), n\ge 0,$ be
the Neumann spectrum of (2.1) with the boundary condition
$y'(0)=y'(1)=0.$ It is well known that $\m_n,\n_n \in
[\a^-_n,\a^+_n]$ and $\n_0\le\a_0^+$. The energy Riemann surface
$\L_E(H)$ for the Hill operator consists of 2 sheets
$\L^{(1)}=\C\sm \cup \s_n $ and $\L^{(2)}$. Each sheet is a copy
of the complex plane slit along  $\s (H)$. The first sheet
$\L^{(1)}$ is glued to the sheet $\L^{(2)}$ by  identifying
"crosswise" all slits in $\s (H)$ (if $p\ev 0$, then we have the
Riemann surface of $\sqrt \l$). We define the quasimomentum
$k(\l)=\arccos\D(\l), \l\in \C_+$, 
which is fixed by $k(0)=0$. The function $\sin k(\l)$
is analytic on the Riemann surface
$\L_E(H)$ (see [M]).  The function $k(\l)$ has an analytic extension from $\C_+$ into  the Riemann surface
$\L_E(H)$ without slits $\g_n(H), n\ge 1$.
Recall that if a gap $|\g _n|=0,$ then
$k(\l)$  is analytic at $\l=\a_n^{\pm}$; if $|\g _n|\neq 0,$ then
$k(\l)$ has branch points $\a_n^{\pm}$ and
$$
k(\l)=\pi n+i\sqrt {-2M_n^{\pm}(\l-\a_n^{\pm})}(1+{\rm o}(1)),\ \ \
{\rm as}\ \ \ \ \l\to \a_n^{\pm}, \ \ \l\in \g _n(H),
$$
where $\pm M_n^{\pm}>0$ is the effective mass; here and below $\sqrt z>0, z>0$ (see [KK2]).
We introduce the Bloch functions $\p_{\pm}(x,\l)$ and the Weyl functions $m^{\pm}(\l)$ by
$$
\p_{\pm}(x,\l)=\vt(x,\l)+m_{\pm}(\l)\vp(x,\l), \ \ \ \ \ \
m^{\pm}(\l)={a(\l)\pm i\sin k(\l)\/\vp(1,\l)}.
$$
We define resonances of $T_t$. There are different kinds of the
resonances. Due to (1.2) $\s_{ac}(T_t)=\s_{ac}(T_0)$ does not
depend on $t\in\R$. It has the decomposition
\[
\s_{ac}(T_t)=\cup_2^4\s^{(n)}, \ \ \s^{(2)}=\s (H_1)\cap\s (H),
\ \ \s^{(3)}=\s (H_1)\sm \s (H),\ \ \s^{(4)}=\s (H)\sm \s (H_1).
\]
In order to describe the resonances we need the energy Riemann
surface $\L_E(T_0)$ for the operator $T_0$. There are 2 cases. In
the first case $\s_{ac}(T_0)\neq \s ^{(2)}$, the Riemann surface
$\L_E(T_0)$ consists of 4 sheets $\L_0^{(1)},
\L_0^{(2)},\L_0^{(3)},\L_0^{(4)}$. Each sheet is a copy of the
complex plane slit along $\s_{ac}(T_0)$. The first sheet
$\L_0^{(1)}$ is glued to the sheet $\L_0^{(n)}$ by identifying the
sides of all slits in $\s^{(n)}, n=2,3,4$. Similarly, the second
sheet $\L_0^{(2)}$ is glued to the sheets $\L_0^{(3)}$ and
$\L_0^{(4)}$ by identifying "crosswise" the sides of all slits in
$\s^{(4)}$ and $\s^{(3)}$ respectively. The third sheet
$\L_0^{(3)}$ is glued to the sheet $\L_0^{(4)}$  along slits in
$\s^{(2)}$. For each complex number $\l\in \L_0^{(1)}$ the number
$\l^{(n)}, n=1,2,3,4$ will denote the corresponding point on the
sheet $\L_0^{(n)}$. In the second case $\s_{ac}(T_0)=\s^{(2)}$ the
energy Riemann surface $\L_E(T_0)$ consists of 2 sheets
$\L_0^{(1)},\L_0^{(2)}$ and coincides with the Riemann surface
$\L_E(H)$ for the Hill operator $H$.

Let $\cB$ denote the class of bounded operators in $L^2(\R)$. For each $\e\in
C_0^{\iy}(\R )$  we introduce the operator-valued function
$A:\L_0^{(1)}\to\cB$ by $A(\l)=\e R(\l )\e $, where $R(\l
)=(T_t-\l )^{-1}$. Assume that $A(\l), \l\in\L_0^{(1)}$, has a
meromorphic continuation across the set $\s^{(m)}$ on "the sheet"
$\L_0^{(m)}, m=2,3,4.$ Suppose that $A(\l)$ has a pole
$\l_r\in\L_0^{(m)}$ (this does not depend on the choice of  $\e$). We call $\l_r$ a resonance of $T_t$. Let $\# ^{(m)}(T_t, \o )$ be the number of resonances of $T_t$ in an interval $\o\ss\L_0^{(m)}, m=2,3,4$.  We formulate our first result.

\no    {\bf Theorem 2.1}.
{\it Let $q_t, t\in \R$, be a biperiodic potential in the sense of (1.1). Then

\no i) Relations (1.2-3) are fulfilled.

\no ii) For each $\e\in C_0^{\iy}(\R )$ the operator-valued
function $ A:\L_0^{(1)}\to\cB$ has a meromorphic continuation into
the Riemann surface $\L_E(T_0)$ described above. The number $\l_e\in \L_0^{(1)}$
is an eigenvalue of the operator $T_t$ iff the same number considered as an element $\l_e^{(2)}\in\L_0^{(2)}$ on
the second sheet  is a resonance of $\wt
T_t$ given by $\wt T_ty =-y''+\wt q_t(x)y$,
where $\wt q_t=\c_-p_2(\cdot+t)+\c_+p_1$.

\no iii) The Riemann surface $\L_E(T_t^{di})$ for the dislocation
operator coincides with the Riemann surface $\L_E(H)$ for  the
Hill operator. 
If $p\not\ev s$, then the Riemann surface $\L_E(T_t^s)$
for the half-solid case consists of 4 sheets.
}

\no    {\it Remark.} 1) To say anything about all complex
resonances on the level of generality of Theorem 2.1 is highly
non-trivial. It is not known, for precisely which subclass of
potentials there exist non real resonances. Theorem 2.1, however,
gives control on all resonances on the second sheet $\L_0^{(2)}$
(but not on the third and fourth), since they are
eigenvalues of the self-adjoint operator $\wt T_t$. In particular,
they are real and belong to $\cup \g_n (T_t)\ss \L_0^{(2)}$. It is
this result which makes the operator $\wt T_t$ important in
Theorem 2.1 (and later on in Lemma 4.1).  2) In the case of the
dislocation operator ($p_1=p_2$)
the energy Riemann surface $\L_E(T_t^d)$ consists of only 2 sheets
$\L_0^{(1)},\L_0^{(2)}$. Thus, by Theorem 2.1, all resonances are
real in this case. 3) In the case of the half-solid ($p_1=const$)
there are (except for $p_1=p_2$) 4 sheets and it can be shown that
in this case there exist complex resonances [KP].

Introduce the subspaces of even potentials $L_{even}^r(\T)\ss L^r(\T)$ given by
$$
L_{even}^r(\T)=\lt\{p\in L^r(\T): p(x)=p(1-x), 0 < x < 1 \rt  \},\ \ \ r\ge 1.
$$
It is well known that a potential $p\in L_{even}^1(\T)$ iff
$|\g_n(H)|=|\m_n(p)-\nu_n(p)|$ for all $n\ge 1$.
Next we consider a biperiodic potential with even potentials $p_1, p_2$.

\no    {\bf Theorem 2.2}.
{\it Let $q_0$ be a biperiodic potential in the sense of (1.1), where $p_1, p\in L_{even}^1(\T)$ and $t=0$.  Assume that some gap of $T_0$
is given by $\wt\g(T_0)=\wt\g_1\cap\wt\g_2\neq\es $
for some gap $\wt\g_n$ of $H_n, n=1,2$.
In the case $|\wt\g_n|=|\wt\m_n-\wt\n_n|<\iy, n=1,2$ we denote by  $\wt\m_n, \wt\n_n\in \ol {\wt\g_n}$
the corresponding Dirichlet and Neumann eigenvalues. Then
\[
\#(T_0, \wt\g(T_0))= \cases{0\ \ \ &if\ \ \ \ $(\wt\m_1-\wt\n_1)(\wt\m_2-\n_2)>0, $\cr
                    1\ \ \ &if\ \ \ \ $(\wt\m_1-\wt\n_1)(\wt\m_2-\wt\n_2)<0. $\cr},\ \
 \ \ \ {\rm if} \ \ \ |\wt\g_1|, |\wt\g_2|<\iy,
\]
%ii) Let $\g _1=(-\iy , 0)$ and let $ \g _2$ be a finite gap. Then
\[
\#(T_0, \wt\g(T_0))=\cases {0,\ \ \ &if\ \ \ \ $\wt\mu _2<\wt\nu_2,$\cr
                      1,\ \ \ &if\ \ \ \ $\wt\n_2<\wt\m _2\le 0. $\cr},
 \ \ \  \ \ \  \ \ \  \ \ \  \ \ \  \ \ \ {\rm if} \ \ \  |\wt\g_2|<\iy=|\wt\g_1|.
\]
Moreover, if $\wt\g (T_0)=(-\iy,0)$, then $\#(T_0, \wt\g (T_0))=0$.}

 We now consider the problem of the reconstruction of $p_1$: if we know $p$ plus some spectral data. In fact we consider the inverse problem  including the problem of characterization for $T_0=-{d^2\/dx^2}+q_0$. We introduce the sequence
$d=\{d_n\}_1^{\iy}$, where $d_n\in \{0,1\}=0$. Define the spaces
$$
Q=\lt\{q_0=\c_-p_1+\c_+p: \ p_1, p\in L_{even}^2(\T),\g_n(H_1)=\g_n(H),
 \ {\rm for\ all} \ n\ge 0\lt\},
$$
$$
P=\lt\{(p,d):\ \ p\in L_{even}^2(\T),  d=\{d_n\}_1^{\iy }, \
d_n\in \{0,1\},\ \ 
\ d_n=0 \ \ {\rm if\ } \ |\g_n(H)|=0 \lt\}.
$$
Define the mapping $\r :Q\to P$ given by $\r (q_0)=(p,d)$,
where $d_n=\#(T_0, \g_n(T_0)), n\ge 1$. Note that (2.3) yields
$d_n\in\{0,1\}, n\ge 1$ and (1.2) implies $d_n=0$ if
$|\g_n(H)|=0$. We shall show that any sequence $d=\{d_n\}_1^{\iy}$
of the number of eigenvalues in the gaps actually occurs. More
precisely, we have the following inverse result:

\no   {\bf Theorem 2.3.}
{\it  The mapping  $\r :Q\to P$ is 1-to-1 and onto.
Moreover, %for $G=(\sum |\g_n(T_0)|^2)^{1/2}>0 $
for $\|p\|^2\ev\int_\T|p(x)|^2dx$  the following estimates are fulfilled:
\[
\|p_1\|=\|p\|\le 2 G(1+G^{1/3}), \ \ G\le 2\|p\|(1+\|p\|^{1/3}),\ \ \
\ \ G=(\sum |\g_n(T_0)|^2)^{1/2}.
\]
Remark.}  i) It follows from assertion of Theorem, that $p_1$
in the decomposition $q_0=\c_-p_1+\c_+p$ is uniquely determined
by $(p,d)$.
ii) Assume that  $q_0\in Q$ and $|\g_n(H)|>0, n\ge 1$ for some $p\in L_{even}^2(\T)$.
Then the definition of $Q$ and (1.2) yield $\g_n(T_0)=\g_n(H)$ for all
$n\ge 1$. If $d_n=1$ for all $n\ge 1$, then Theorem 2.3 yields
$\#(T_0, \g _n(T_0))=1,
n\ge 1,$ and $\#(T_0, \R)=\iy $. If $d_n=0, n\ge 1$, then Theorem 2.3
yields $\#(T_0, \g _n(T_0))=0, n\ge 1,$ and $T_0=H$.

We now consider the {\bf dislocation operator}
$$
T_t^{di}=-{d^2\/dx^2}+p_{(t)}(x)\ \ {\rm in}\ \ L^2(\R ),
\ \ \ p_{(t)}=\c_-p+\c_+p(\cdot+t), \ \ p\in L^1(\T), \ t\in \R.
$$
 Using (1.2-3) we obtain the following relations
\[
\s (T_t^{di})=\s_{ac}(T_t^{di})\cup \s_{d}(T_t^{di}),\ \ \  \s_{ac}(T_t^{di})=\s (H),
\ \ \ \s_d(T_t^{di})\ss\cup\g_n(H), \ \ \ \# (T_t^{di},\g_n(T_t^{di}))\le 2.
\]
We emphasize that $\g_n(T_t^{di})=\g _n(H)$ for any $t\in [0,1], n\ge 0$, i.e., $\g _n(T_t^{di})$ are independent of $t$. This follows from Theorem 2.1.i.
Then $\s_{ac}(T_t^{di})=\s^{(2)},$ and $\s^{(3)}=\s^{(4)}=\es $, and
the Riemann surface $\L_E(T_t^{di})$ coincides with the Riemann
surface $\L_E(H)$, and consists of 2 sheets
$\L_0^{(1)},\L_0^{(2)}$. 
Introduce the function (here and below $\dot u={\pa\/\pa t}u$)
\[
L(t,\l)=\mp [\dot \P(t,\l)^2-(p(t)-\l)\P(t,\l)^2],\ \ \
\l=\a_n^{\pm},\ \  \ n\ge 0,\ \ t\in [0, 1].
\]
Let $W^2_r(a,b)$ be the Sobolev space of functions $f$
on the interval such that $f^{(r)}\in L^2(a,b), r\ge 0$.
We have the following result about eigenvalues of the dislocation operator.

\no   {\bf Theorem 2.4.}
{\it Let the dislocation potential $p_{(t)}=\c_-p+\c_+p(\cdot+t)$, where $p\in L^2(\T), t\in\R$. Then for each gap $\g_n(H)=(\a_n^-,\a_n^+)\neq\es, n\ge 1$
there exists a function $z_n^{\pm}(\cdot )\in W^2_1(-\ve_n, \ve_n),$  for some
$\ve_n>0$, such that: $z_n^{\pm}(0)=0$ and $\l_n^{\pm}(t)\ev\a_n^{\pm }\mp z_n^{\pm}(t)^2\in \g_n(H)=\g_n(T_t^{di})$ satisfies

\no i) If $z_n^{\pm}(t)>0$ for some $t\in (-\ve_n,\ve_n)$, then
$\l_n^{\pm}(t)\in \g_n(H)$ is an eigenvalue of $T_t^{di}$.

\no If $z_n^{\pm}(t)<0$ for some $t\in (-\ve_n,\ve_n)$, then
$\l_n^{\pm}(t)\in\g_n(H)\ss\L_0^{(2)}$ is a resonance of $T_t^{di}$.

\no Moreover, the following  asymptotics are fulfilled:
\[
z_n^{\pm}(t)=\sqrt{|M_n^{\pm}|/2}\int _0^tL(s, \a_n^{\pm})ds
+O(t^{3/2}),\ \ \ \ \ t\to 0,
\]
\[
z_n^{\pm}(t)=\mp t\sqrt{|M_n^{\pm}/2|}\dot \P(0, \a_n^{\pm})^2+O(t^2),
\ \ \ \ \ t\to 0,\ \ \  \ \ \ {\rm if} \ \ \ \m_n(p)=\a_n^{\pm}.
\]
ii) If, in addition, $p$ is real analytic, then $z_n^{\pm}(\cdot )$ is real analytic on $ (-\ve_n, \ve_n)$.

\no Remark.} i) If in (2.9) we have  $\m_n(p)=\a_n^-,$ then
$z_n^{-}(t)>0$ for $t>0$ and $\l _n^{-}(t)$ is an eigenvalue;
if $\m_n(p)=\a_n^+,$ then $z_n^{-}(t)<0$ for $t>0$ and $\l_n^{-}(t)$ is a resonance.
ii) Let $p$ be smooth potential  and $n\gg 1$.
 Then using (2.8) we deduce  that  $z_n^{-}(t)>0,$ as $ t>0$ (we get an eigenvalue) and  $z_n^{+}(t) <0$ (we get a resonance). Hence for fixed $n\gg 1$ and any small $t>0$ we have an eigenvalue near $\a_n^-$ and a resonance near $\a_n^+$ on the second sheet.
iii) It is possible to formulate  Theorem 2.4 for potentials
 $p\in L^1(\T)$ and then the asymptotics in (2.8) has the form
$z_n^{\pm}(t)= \sqrt{|M_n^{\pm}/2|}\int _0^tL(s, \a_n^{\pm})ds +O(t),$ as $t\to 0$

Using Theorem 2.4 we construct the dislocation
operator $T_t^{di}$, such that $\#(T_t^{di}, \g _n(T_t^{di}))=2$ for all $n=1, .., N$ for each fixed $N\ge 1$.

\no   {\bf Theorem 2.5.}
{\it  For any finite sequences $d=\{d_n\}_1^{N}, \{s_n\}_1^{N}, d_n\in\{0,1\}, s_n>0,  N\ge 1$   there exists a potential $p\in L_{even}^2(\T)$ and $\ve >0$ such that the each dislocation operator $T_t^{di}, t\in (0, \ve )$
has  gaps with lengths  $|\g _n(T_t^{di})|=s_n,$ and
$\#(T_t^{di}, \g _n(T_t^{di}))=2d_n$ for any $n=1, 2, .., N$.}

We now consider a {\bf half-solid}, that is the Schr\"odinger operator
$$
T_t^s=-{d^2\/dx^2}+q_t^s(x) \ \ \ {\rm in}\ \ L^2(\R ),
\ \ \ \ q_t^s(x)=s\c_-(x)+\chi_+(x)p(x+t),
$$
where $s,t\in \R$ and $p\in L^1(\T )$ is real. Theorem 2.1
yields :
\[
\s (T_t^s)=\s_{ac}(T_t^s)\cup\s_{d} (T_t^s),\ \ \ \ \ \ \
\s_{ac}(T_t^s)=\s (H)\cup [s,\iy ),
\]
\[
\g_n(T_t^s)=\g_n(H)\cap (-\iy,s), \ \ \ \  \ \#(T_t^s,\g_n(T_t^s))\le 2,
\ \ \ n\ge 0,\ t\in[0,1].
\]
Hence Theorem 2.1 yields if $p\not\ev s$, then the Riemann surface $\L_E(T_t^s)$ consists of 4 sheets.
Note that if $s\le\a_1^-(H)$, then $\s_{ac}(T_t^s)=(s_-,\iy )$,
where $s_-=\min (0,s)$. If $s> \a_1^-(H)$, then there exists a gap
in the spectrum of $T_t^s.$ Our goal is to study the eigenvalues in the gaps $\g _n(T_t^s), n\ge 0,$ and to find how these eigenvalues depend on $t, s$.  It is clear that they depend on $t$ periodically. Define sets $Q_N=\lt\{q_0^s=s\c_-+\c_+p_2:$  \ \ \ $p\in L_{even}^2(\T),\ \ \a_{N}^+<s<\a_{N+1}^-\lt\}$ and
$P_N=\lt\{(r,d,\ve): \{r_n\}_1^{\iy}\in \ell^2, d=\{d_n\}_1^{N}$,
where for any  $n\le N$ the number $r_n\ge 0$ and $d_n=0$ if $r_n=0,
 \ve\in (0,1)\lt\}$  for some integer $N\ge 0$ .
Define the mapping $\o :Q_N\to P_N$ by $\o (q_0^s)=(r,d,\ve )$, where
$$
r_n=|\g_n(T_0^s)|, \ d_n=\#(T_0^s, \g_n(T_0^s)), \
n\le N, \ \  r_n=|\s_n^1(T_0^s)|\sign (\m_n(p)-\n_n(p)), \ \ \ n>N,
$$
$
\ve={\a_{N+1}^--s\/\a_{N+1}^--\a_{N}^+}$ and  $\s_n^1(T_0^s)=(s,\iy)\cap\g_n(H)$
is a segment of the spectrum of $T_0^s$ with multiplicity one.
We consider a half-solid with an even potential $p$.

\no    {\bf Theorem 2.6}
{\it i) Let $T_0^s=-{d^2\/dx^2}+q_0^s(x)$ acting in $ L^2(\R )$, where $q_0^s=s\c_-+\c_+p$ and $p\in L_{even}^2(\T)$.
Then $\#(T_0^s,\g_0(T_0^s))=0$. Let, in addition, $s>a_{m}^+(p)$ for some $m\ge 0$. Then $\s^{(4)}\neq\es$, and for each $\g_n(T_0^s)\neq \es, n=1, 2,..,m,$
the following identities are fulfilled:
\[
\#^{(4)}(T_0^s,\g_n(T_0^s))+\# (T_0^s,\g_n(T_0^s))=1,\ \ \
\#(T_0^s,\g_n(T_0^s))=\cases{1\ \ \ &if\ \ \ \ $\m_n>\n_n,$\cr
                               0\ \ \ &if\ \ \ \ $\n_n>\m_n.$\cr}
\]
ii) Each mapping  $\o :Q_N\to P_N, N\ge 0,$ is 1-to-1 and onto.}

\no In this Theorem we proved that for any $\{r_n\}_1^{\iy}\in \ell^2, d=\{d_n\}_1^{N}$, and $\ve \in (0,1)$ (such that $r_n\ge 0$ and $d_n=0$
if $r_n=0,1\le n\le N$), there exists a unique half-solid potential $q_0^s=s\c_-+\c_+p$ with $p\in L_{even}^2(\T)$.

For fixed $\g_n(T_t^s)\neq \es, n\ge 0,\ s>\a_n^{\pm}$ we introduce the equation
\[
\P_y(y,\a_n^{\pm})=\sqrt{s-\a_n^{\pm}}\P(y,\a_n^{\pm}),\ \ \ \ \ \ y\in [0, 1].
\]
For each  $n\ge 1$ there exist $N\ge n$ roots $y_1,..,y_N\in [0,1]$ of Eq. (2.13) on the interval $[0,1]$, see below Lemma 3.4.
Let $m^{\pm}(\l, t)$ be the Weyl function for the potential $p(x+t)$
and let $\a_0^-\ev-\iy$.

\no    {\bf Theorem 2.7.}
{\it Let $T_t^sy=-y''+q_t^sy$ be an operator acting in $L^2(\R )$, where $q_t^s=s\c_-+\c_+p(\cdot+t)$ and $p\in L^1(\T),t\in \R$. Then

\no i) Suppose that $\a_n^-<\a_n^+<s$ for some  $n\ge 0$ and $y\in [0,1]$ is some root of (2.13).
Then there exists a unique function $z_n^{\pm}(\cdot )\in C(-\ve, \ve ), z_n^{\pm}(0)=0$ for some $\ve >0$, such that :

\no if $z_n^{\pm}(t)>0$ for some $t\in (-\ve,\ve )$, then
$\l_n^{\pm}(t)\ev\a_n^{\pm }\mp z_n^{\pm}(t)^2\in \g_n$
is an eigenvalue of $T_{y+t}^s,$

\no if $z_n^{\pm}(t)<0$ for some $t\in (-\ve, \ve )$, then
$\l_n^{\pm}(t)\in \g_n \ss\L_0^{(4)}$ is a resonance of $T_{y+t}^s$.

\no Moreover, the following asymptotics is fulfilled:
\[
z_n^{\pm}(t)=\sqrt{|M_n^{\pm}/2|}\P (y,\a_n^\pm)^2\int_0^t[p(y+\t)-\a_n^{\pm}]d\t+O(t),
 \ \ \ t\to 0.
\]
\no ii) Suppose that $\a_n^-<s<\a_n^+,$ for some $ n\ge 1,$ or $\nu_0^-<s\le\a_0^+.$
Let $y$ be  a root of the equation $m^+(s, y)=0, 0\le y\le 1$.
Then there exists a unique function $z(\cdot )\in W^2_1(-\ve, \ve ), z(0)=0$
for some $\ve >0$ such that  

\no if $z(t)>0$ for some $t\in (-\ve, \ve )$ then
$\l(t)\ev s+z(t)^2\in \g_n$ is an eigenvalue of $T_{y+t}^s,$

\no if $z(t)<0$ for some $t\in (-\ve,\ve )$ then
$\l(t)\in \g_n \ss\L_0^{(3)}$ is a resonance of $T_{y+t}^s$.

\no Moreover, the following asymptotic estimates are fulfilled:}
\[
z(t)=\int _0^t[p(y+\t)-s]d\t+O(t),\ \ \ \ \ t\to 0.
\]

\section{ Hill operator and translations}
\setcounter{equation}{0}

Let $\vp(x,\l,t), \vt(x,\l,t)$ be the solutions of the equation
\[
-y''+p(x+t)y=\l y, \ \ \ \l\in \C ,\ \ \ \  t\in \R ,
\]
satisfying   $\vp_x(0,\l,t)=\vt (0,\l,t)=1,$ and $\vp(0,\l,t)=\vt_x(0,\l,t)=0$.
Remark that the Lyapunov function $\D (\l)$ for (3.1) coincides with
the Lyapunov function for (2.1) (see [L]). Then  we define the quasimomentun $k$ for (3.1) and again $k(\l)$ does not depend on $t$. Let $\m_n(p,t), n\ge 1, $ be the Dirichlet spectrum of $p(x+t),$ i.e
the spectrum of (3.1) with the boundary condition $y(0)=y(1)=0$
and let $\n_n(p,t), n\ge 0,$ be the Neumann spectrum of $p(x+t)$,  that is
the spectrum of (3.1) with the boundary condition $y'(0)=y'(1)=0.$
We need some results on Eq. (3.1)  (see [L], [PTr], [T],[K2] ).
It is well known that $\vp(1, \l, t)$ is an entire function of $\l$ of order ${1 \/ 2}$ for fixed  $t$. The zeros of  $\vp(1, \l, t)$
coincide with the Dirichlet eigenvalues $\m_n(p,t), n\ge 1,$ and
the following asymptotics are fulfilled:
\[
\m_n(p,t) =(\pi n)^2+\int _0^1 p(x)dx-\int_0^1p(x+t)\cos 2\pi nxdx +O(1/n),
\ \ \ \ \ \ {\rm as} \ \ \ \ n\to \iy ,
\]
uniformly on bounded subsets of $[0,1]\ts L^1(0,1)$.
It is well known that $\vt_x(1,\l,t)$ is an entire function of $\l$
of order $1/2$ for fixed  $t$. The zeros of  $\vt_x(1,\l,t)$
coincide with the Neumann eigenvalues $\n_n(p,t), n\ge 0.$
The functions $\m_n(p,t), \n_n(p,t)$  are 1-periodic. If the parameter $t$ runs through the interval $[0,1]$, then $\m_n(p,t), \n_n(p,t)$ run through
the gap $\g_n=(\a_n^-,\a_n^+), n\ge 1$. If the gap $\g_n=\es$, then
$\m_n(p,t), \n_n(p,t)$ don't move and $\mu _n(p,t)=\nu_n(p,t)=\a_n^{\pm}$.
The eigenvalue $\nu_0(p,t)\le \a_0^+, t\in [0, 1].$
In the book [PTr] there are the asymptotics of the solutions
$\vp(x,\l), \vt(x,\l)$ as $|\l|\to \iy .$ Repeating
it for Eq. (3.1) we obtain the asymptotics for
$\vp(x,\l,t), \vt(x,\l,t)$ as $|\l|\to \iy$. For example,
\[
\vp(x,\l,t)={\sin \sqrt{\l}x\/\sqrt{\l}}+O({\exp{\Im \sqrt{\l}x} \/ \l}) \ \ \ \ \ {\rm as}\ \ \ \ x, t\in [0, 1],\ \  |\l|\to \iy .
\]
These asymptotics can be differentiated with respect to $x, t $
and /or $\l$ and are uniform on $[0, 1]\ts [0, 1]\ts L^1(0, 1)$.
We have the Trubowitz identity (see [Tr])
\[
\vp(1,\l,t)=\vp(1,\l)\p_+(t,\l)\p_-(t,\l)=\prod_{n\ge 1}{(\m_n(p,t)-\l)\/(\pi n)^2},
\ \ \ \l\in \C , \ \ t\in \R,
\]
and the equality $\vp(1, \l)\p_+(t,\l)\p_-(t,\l)=-2\D'(\l)\P (t,\l)^2$
at $\l=\a_n^{\pm}$, yields
\[
\vp(1,\a_n^{\pm},t)=(-1)^n2M_n^{\pm}\P (t,\a_n^{\pm})^2, \ \ t\in \R, \ \ n\ge 0.
\]
where the eigenfunction $\P (t,\a_n^{\pm})$ is defined at the start of Section 2. Define a function $\f_n(\l,t)$ by $\vp(1,\l,t)\ev (-1)^n(\l-\m_n(p,t))\f_n(\l,t)$. Then we have
\[
\f_n(\l, t)>0, \ \ \ \l\in \g _n,   \ \ {\rm and}\ \
(-1)^n\vp(1, \l, t)>0,\ \ \  \l>\mu _n(t),\ \  \l\in \g _n,
\]
Introduce a function $b(\l)=-i\sin k(\l)$.
It is known that 
\[
k(\l)\ev \pi +iv(\l),\ \  v(\l)>0, \ \ \  b(\l)=(-1)^n\sqrt{\D^2(\l )-1}, \ \
\ \ \ \l\in \g _n \in \L ,
 \]
where $\sqrt{\D (\l)^2-1}>0$ as $\l\in \g _n\ss \L$.
Below we need the identity (see [KK1])
 \[
M_n^{\pm}=-\D(\a_n^{\pm})\D'(\a_n^{\pm}), \ \ \ n\ge 0,
 \]
and using (3.7-8) we obtain
\[
b(\l)=(-1)^n z( \sqrt{|2M_n^{\pm}|}+O(z^2)) ,\ \ {\rm as} \ \
\l\in \g _n\neq \emptyset ,\ \ \l=\a_n^{\pm}\mp z^2, \ \ z\to 0.
 \]
The Weyl function $m^{\pm}$ for the potential $p(x+t)$ has the form
\[
m^{\pm}(\l,t)={a(\l,t){\pm}i\sin k(\l)\/\vp(1,\l,t)}
={a(\l,t){\mp}b(\l)\/\vp(1,\l,t)}, \ \ \ \ a(\l,t)\ev {\vp_x(1,\l,t)-\vt(1,\l,t)\/2}.
\]
Below we need the identity
\[
a^2(\l, t)+1-\D^2(\l)=a^2(\l,t)-b^2(\l)=-\vp(1,\l,t)\vt_x(1,\l, t).
\]
Recall that $\dot u={\pa\/\pa t}u$. We have the equations
\[
\vcenter{
\hbox{$\dot\vt_x(1,\l,t)=(\l-p(t))\dot\vp(1,\l,t),\ \ \ \dot\vp(1,\l,t)=2a(\l,t),$}
\vskip\baselineskip
\hbox{$\kern2.8em\dot a(\l,t)=-\vt_x(1,\l,t)-(\l-p(t))\vp(1,\l,t),$}}
\]
(see [L]) and the following identities
\[
\vp(x, \l, t)=\vt(t, \l) \vp(x+t, \l)-\vt(x+t, \l) \vp(t, \l),
\]
\[
\vt (x, \l, t)=\vt(x+t, \l) \vp_x(t, \l)-\vt_x(t, \l) \vp(x+t, \l).
\]
Let $C^m(\T), m\ge 0,$ be the space of m times continuously
differentiable real-valued 1-periodic functions. Suppose that
$p\in C^1(\T),$  then for any $t\in [0,1]$ the identity (the trace
formula)
\[
p(t)=\a_0^++ \sum _{n\ge 1}\lt(\a_n^-+\a_n^+-2\mu _n(p,t)\lt),\ \
\]
holds, where the series converges absolutely and uniformly (see [L]).
We need the following result on the Dirichlet spectrum (see [K2]).

\no    {\bf Theorem 3.1.}
{\it  Let a real potential $p \in L^1(\T)$ and $\m_n(t)=\m_n(p,t), n\ge 1$. Then

\no  i)  Each $\mu _n(\cdot )\in C^2(\T), n\ge 1,$  and
$\mu _n''' \in L^1(\T).$
Let in addition $p\in L^2(\T)$ (or $p\in C^m(\T), m\ge 0$).
Then  $\mu _n'''\in L^2(\T)$ (or $\mu _n(\cdot )\in C^{m+3}(\T), m\ge 0$).

\no  ii) There exists a function $y_n\in C^1(\R ),$ such that
$\m_n(t)=\a_n^-+|\g_n|\sin^2 y_n(t),$ where
$\mu _n(0)=\a_n^-+|\g_n|\sin^2 y_n(0),$ and uniformly on $t\in [0, 1]$
the  following asymptotics are fulfilled:
\[
 y_n(t)=y_n(0)+\pi n t+O(1/n), \ \ \ \dot y_n(t)=\pi n +o(1),\ \ \
 \ \ \ \ n\to \iy .
\]
\no iii) Suppose that $\mu _n(t_0)=\a_n^-$ or $\mu _n(t_0)=\a_n^+,$ for some $t_0\in [0, 1], $ and $ n\ge 1.$ Then the following asymptotics is fulfilled:
\[
\mu _n(t_0+t)=\mu _n(t_0)+t^2\ddot \mu _n(t_0)/2+o(t^2),\ \ \
   \ddot \mu _n(t_0) =-4M_n^{\pm}/ \vp'_\l(\mu _n(t_0), t_0)^2.
\]
\it Remark.} i) In other words we have the following result.
Slit the n-th gap $\g _n\neq \es$ and place $\mu _n$ on the upper or
lower lip according to the signature of $\sinh q(\mu _n)$, i.e., on the upper when positive and on the lower when negative. Then $\mu _n(t)$ runs clockwise
around the "circle ", changing lips when it hits $\a_n^{\pm}$, making $n$
complete revolutions in unit time. Then (roughly speaking)
$\m_n(t)=\a_n^-+|\g_n|\sin^2 \pi nt,$ when $t$ runs through the interval $[0,1].$

In order to study the eigenvalues in gaps we need properties
of the function $\z(\l,t)=\dot\vp(1,\l,t)/(2\vp(1,\l,t))$.

\no    {\bf Lemma 3.2.}
{\it For each $(t,p,\l)\in F= [0,1]\ts L^1(\T)\ts \C\sm \{\m _n(p,t),n\ge 1\}$
the following  identities are fulfilled:
\[
\z (\l,t)\ev{\dot\vp(1,\l,t)\/2\vp(1,\l,t)}=
{1\/2}\sum_{n\ge 1}{\dot\m_n(p,t)\/\m_n(p,t)-\l},\ \ \ \
\]
\[
\dot m^{\pm}(\l,t)=(p(t)-\l)-m^{\pm}(\l,t)^2,
\]
\[
\dot\z (\l,t)=(p(t)-\l)-\z^2(\l,t)-{b^2(\l)\/\vp^2(1,\l,t)},
\]
\[
\int_0^1\z (\l,t)dt=0,\ \ \ \l<\a_1^-,\ \ \ \
\int_0^1\z (\a_0^+,t)^2dt=\int_0^1p(t)dt-\a_0^+,
\]
\[
\vp(1,\l,t)\dot\z (\l,t)=\pm (-1)^n2M_n^{\pm}L(t,\l),\ \ if \ \l =\a_n^{\pm}\neq\m_n(p,t).
\]
where the series converges absolutely and uniformly on compact sets in $F$ and $L$ is given by (2.7).

\no Proof.}  Let $\m_n(t)=\m_n(p,t)$. The functions $\vp(1,\l,t)$ and $\dot\vp(1,\l,t)$
are entire in $\l$. xThe zeros of
$\vp(1, \l, t)$ have the asymptotics (3.2). Using (3.2-3) we get 
the asymptotics  $\z (\l, t)=O(1/\sqrt{\l})$ as $|\l|\to \iy ,
 |\sqrt{\l}-\pi n|\ge 1/4$. Then (3.4)  yields 
$$
\z (\l,t)-\sum_{n\ge 1}^N{\dot\m_n(t)\/\m_n(t)-\l}=
{1\/2\pi i}\int_{|\sqrt{z}|=\pi (2N+1)/2}{\z (z,t)dz\/z-\l},
$$
since the residue of the function $\z (z,t)$ at the simple pole $\m_n(t)$ has the form: ${\rm Res}\z (\l,t)=\dot\vp(1,\m_n(t),t)/\vp'(1,\m_n(t),t)=\dot\m_n(t)$.
Then as $N\to \iy $ we get (3.18). Remark that (3.16) yields
$\dot \mu _n(t)=O(n|\g_n|)$ as $n\to \iy $ uniformly on $t\in [0, 1].$ Then by (3.2), series (3.18) converges absolutely and uniformly on compact sets.

We will obtain the equations for the function $m^{+}(\l,t)$,
the proof for $m^-(\l, t)$ and (3.20) is similar.
Using (3.12), (3.10-12) we deduce that
$$
\dot m^+={\dot a\/\vp}-{(a-b)\dot\vp \/\vp^2}=(p-\l)-{\vt_x\/\vp}-{2(a-b)a\/\vp^2}
=(p-\l)+{a^2-b^2\/\vp^2}-{2am^+\/\vp},
$$
which yields (3.19). Integrating $\z=\dot\vp/2\vp$ and since
$\vp(\l, t)>0$ if $\l<\a_1^-, t\in [0, 1],$ we obtain the first
identity in (3.21). Integrating (3.20) implies the second
identity in (3.21). Substituting (3.5) into (3.20) we  have (3.22). $\BBox$

Now we will obtain more exact estimates concerning $L(t,\a^{\pm}_n)$ defined
in (2.7). These results are used to prove the existence of two eigenvalues in gaps (see Theorem 2.5).  Recall that if $p\in L_{even}^2(0,1)$, then $\n_0<\m_1$ and $\g_n(H)=(\m_n,\n_n)$ or $\g_n(H)=(\n_n,\m_n), n\ge 1$ [GT].

\no    {\bf Lemma 3.3.}
\no {\it i) Let $p\in C^2(\T)$ be even. Assume that for
some $N\ge 1$ and for all $n>N$ the gap $\g_n(H)$
are given by $\g_n(H)=(\m_n(p),\n_n(p))$. Then
\[
p(0)-\a_N^{\pm}\ge -(\pi N)^2+\sum_{m\ge N+1}|\g _m|
-2\sum _1^{N}|\g_m|.
\]
ii) For any finite sequences $\{s_n\}_1^N, \{d_n\}_1^N,$ where $s_n>0, d_n\in\{0,1\}$, there exists $p\in L_{even}^2(\T)$ 
with gap lengths $|\g_n(H)|=s_n$
and $(-1)^{d_n}L (0,\a_n^{\pm})<0$ for all $n=1, \dots , N.$

\no Remark.} Roughly speaking  the result of ii) is the effect of a big gap $\g_{N+1}$ such that $|\g_{N+1}|>(\pi N)^2+2(|\g_1|+...+|\g_{N}|)$. It is important that we have 2 types
of gaps: 1)  the gap  $\g_n(H)=(\m_n(p),\n_n(p))$ for all $n>N$, 2) 
$\g_n(H)=(\m_n,\n_n)$ or $\g_n(H)=(\n_n,\m_n)$, which depends on
$d_n$ for $1\le n\le N$.

\no {\it Proof.}
i) Using the trace formula (3.14) and the identity $\a_n^+=
\sum _{m\ge 1}^{n}(|\s_m|+|\g_m|)$ we obtain
$$
p(0)-\a_N^+=\sum_{m\ge 1}(\a_m^++\a_m^--2\mu_m(p))-\sum_{m=1}^{N}(|\s_m|+|\g_m|)
\ge \sum_{m\ge N+1}|\g_m|-2\sum_{m\ge 1}^{N}|\g_m|-\sum_{m=1}^{n}|\s_m|
$$
and the estimate $|\s_m|<\pi ^2 (2m-1)$ (see [Mos]) yields (3.23).

\no ii) We need a result from [GT]. For an even periodic potential we define a signed
 gap length $l_n=\mu _n-\nu_n, n\ge 1, $ and the corresponding sequence
$l= \{l_n\}^{\iy }_1$. For any sequence
$\{t_n\}^{\iy }_1\in \ell^2$ there
exists a unique even periodic potential $p$  such
that the signed gap length $l_n=t_n, n\ge 1. $
Using this result we fix the number $N$ and due to (3.23) we take
a sequence of signed gap lengths $\{l_n\}^{\iy }_1$
such that  $p(0)-\a_N^+>0$. Hence  for all $n=1,  .., N$ we get
$p(0)-\a_n^{\pm}>0$. Thus for each $1\le n\le N$ we obtain: 

If $d_n=1$ we take $\mu _n=\a_n^-$ and
$\nu_n=\a_n^+$, then $L (0,\a_n^{\pm})>0$.

If $d_n=0$ we take $\n_n=\a_n^-$ and $\m_n=\a_n^+$, then $L(0,\a_n^{\pm})<0.
\ \ \BBox$

We need some results about the roots of the equation $m^{\pm} (\l, t)= \omega ,t\in [0, 1],$ and some formulas for even potentials.

\no    {\bf Lemma 3.4.}
{\it i)  Let $p\in L^1(\T)$ and a gap $\g _n(H)\neq \es$ for some $n\ge 1$. Then for any fixed  $(\l,\o)\in [\a_n^-,\a_n^+]\ts\R$, there exist $N_{\pm}\ge n$ roots of the equation $m^{\pm} (\l, t)=\omega ,t\in [0, 1].$

\no ii) For any even potential $p\in L^1(\T)$ the following
identities are fulfilled:
\[
a(\l, 0)\ev 0,\ \ \ \  \D(\l)\ev \vt(1, \l),\ \   \ \
m^{\pm}(\l)\ev {\pm i\sin k(\l) \/ \vp(1, \l ,0)}.
\]
\no  Proof.} Let $\m_n(t)=\m_n(p,t), n\ge 1$.
i) We consider $m^+$ and firstly let $\l=\a_n^-,$ the proof for $m^+$ and
$\l=\a_n^+$  is similar. By Theorem 3.1,  $\m_n(\cdot)\in C^2(\T)$
and there exist points $\t_r\in [0,1), r=1,.., n,$ such that $\m_n(\t_r)=\a_n^-$ and
$\dot\m_n (\t_r)=0, \ddot\m_n (\t_r)>0.$ By (3.18),
$$
m^+(\l,t)=\z(\l,t)={\dot\m_n(t)\/\m_n(t)-\a_n^-}+\sum_{m\neq n}
{\dot\m_m(t)\/\m_m(t)-\a_n^-}.
$$
Then $m^{+}(\l,\cdot )$ maps the interval $I_r=(\t_r,\t_{r+1})$
onto the real line $\R $. Therefore, for any number $\o \in \R $
there exist $N\ge n$ roots of the equation $m^+(\l, t)=\o $.

Secondly, let $\l\in (\a_n^-,\a_n^+)$ and $b(\l)>0,$ the proof for $b(\l)<0$
is similar. We get $2a(\l,t)=\dot \vp (1,\l,t)=\dot \m_n(t)(-1)^n\f(\l,t)$
at $\l=\m_n(t).$ Since $\m_n(t)$ crosses the point $\l$ exactly $2n$ times
there exist points $\t_r\in [0,1)$ $r =1,2,3, .., n,$ such that $\m_n(\t_r)=\l$ and
$(-1)^n\dot\m_n (\t_r)<0.$ Then $a(\l,\t_r)-b(\l)\neq 0$, and we have
\[
m^+(\l,t)={a(\l,t)-b(\l)\/\vp(1,\l,t)}={a(\l,t)-b(\l)\/(-1)^n(\l-\m_n(t))\f_n(\l,t)}.
\]
Recall that each function $\f_n(\l,\cdot ), \l \in \g_n$ has no zero on
the interval $[0,1]$. By the properties of $\m_n(t)$ from Theorem 3.1, the function
$m^{+}(\l,\cdot )$ maps the interval $I_p=(\t_p, \t_{p+1})$
onto the real line $\R $ .  Then for any $\o \in \R $ there
exist $N_{\pm}\ge n$ roots of  the equation $m^+(\l, t)=\o$.

\no ii) It is  well known that for an even
potentials the points $\mu _n(0)$ and $\nu_n(0)$ lie on the endpoints of the gap $\g _n$(see [GT]) . Then by Theorem 3.1, $\dot \mu _n(0)=0$, and using (3.18)  we get
$ \z (\l , 0)=0$ for all $  \l\in \C $. Hence relations (3.10), (3.12) imply
(3.24).  $\BBox$

We consider the equation $\z (0, t)=r, t\in [0, 1),$ for fixed
$r\in \R.$ Let $\p (t)$ be any smooth 1-periodic positive
function. Then we take $\P(t, 0)=\p(t)$ and the potential
$p(t)=\ddot \p(t)/\p(t).$ In this case the function $\z (0,t)$ is
periodic smooth with $\int _0^1\z (0,t)dt=0$ see (3.21). We are
able to get a function $\p$, when the equation  $\z (0, t)=r$ has
any number of roots (depending on $\p$). In order to study the
ground state we need some properties of the Weyl function.

\no    {\bf Lemma 3.5.}
{\it Let $p\in L^1(0, 1).$ Then $\pm m^{\pm}(\l)\to -\iy$
as $\l\to -\iy ,$ and
\[
\pm m^{\pm}(\l)<0 \ \ \ {\rm if}\ \ \l<\n _0,\ \ \ {\rm and}\ \
m^{+}(\l)m^{-}(\l)>0, \ \ \ {\rm if} \ \ \  \n _0<\l< 0,
\]
Moreover, if $\pm m^{+}(0)>0$, then $\pm m^{+}(\l)>0$ for $\n _0<\l< 0$.

\no   Proof.} (4.5) yields the asymptotics $\pm m^{\pm}(\l)\to -\iy $ as $\l\to -\iy$. We have the
identity  $a^2-b^2=-\vp\vt_x$. It is well known that: $\vp(\l)>0,
\l<\m_1,$ and $\vt_x(\l)>0, \l<\n_0$ and $\vt_x(\l)<0, \n_0<\l<0$.
Hence we get (3.26) since $m^{\pm}$ has no zero on the intervals
$(-\iy,\n_0)$ and $(\n_0,0)$.

Let $\pm m^{+}(0)>0$. Due to (3.11), the function $m^+$ has no zero on the interval $(\nu_0, 0)$.
Then for the case $\l\in (\nu_0, 0)$ we have $ m^+(\l)>0$, if $ m^+(0)>0$, and $m^+(\l)<0,$ if $m^+(0)<0.
\ \ \BBox $

\section{Biperiodic potentials}
\setcounter{equation}{0}

In this section we study the spectrum of $T_t$ at $t=0$, i.e.,
$T_0=-{d^2\/dx^2}+q_0$, where $q_0=p_1\c_-+\c_+p$
and $p\in L^1(\T),p_1\in L^1(\t\T)$.
 We now determine the eigenfunctions of $T_0$. Let $u(x,\l),
v(x,\l)$ be the solutions of the equation
\[
-y''+q_0y=\l y, \ \ \ \l \in \C,
\]
satisfying   $u'(0,\l)=v(0,\l)=1,$ and $u(0,\l)=v'(0,\l)=0.$ Let $\vp_j(x, \l),
\vt_j(x, \l), j=1, 2,$ be the solutions of the equation
$$
-y''+p_jy=\l y, \ \ \ \l\in \C ,
$$
satisfying $\vp _j'(0,\l)= \vt_j(0,\l)=1,$ and $
\vp_j(0,\l)=\vt_j'(0,\l)=0.$ Remark that $\vt , \vp ,\vt _j,
\vp_j$ are entire functions of $\l\in \C ,$ and are real on the
real line. Then we obtain
\[
 v(x,\l)=\left\{\begin{array}{cc} \vt_1(x, \l)\ \ {\rm if}\ \  x < 0,\\
                                  \vt_2(x, \l)\ \ {\rm if}\ \ x > 0,
\end{array}\right.
\ \ \ \ \ \
 u(x,\l)=\left\{\begin{array}{cc} \vp_1(x,\l)\ \ {\rm if}\ \ x < 0,\\
                                  \vp_2(x,\l)\ \ {\rm if}\ \ x > 0,
\end{array}\right.
\]
Let $\D_j(\l)$ be the Lyapunov function for the potential $p_j,
j=1,2$. We need also the quasimomentum $k_j(\l)= \arccos
\D_j(\l)$, which is analytic in the domain $\L_j=\C \sm \s (H_j)$.
Introduce the Bloch functions $\p^{\pm}_{j}(\cdot,\l)\in
L^2(\R_{\pm}), \l\in\L_j,$ for the operator $H_j=-{d^2\/dx^2}+p_j,
j=1, 2$, by
\[
\p ^{\pm}_{j}(x, \l) =\vt _j(x, \l)+ m_j^{\pm }(\l)\vp _j(x, \l),
\ \ \ \ m_j^{\pm}(\l)={a_j(\l)\pm i\sin k_j(\l) \/ \vp _j(\l)}.
\]
Then for $\l\in\L_j, j=1,2$ we have
\[
\vp_j(x, \l)={\p _j^+(x, \l) -\p ^-_j(x, \l)\/ w_j(\l)},
\ \ \vt_j(x, \l)={m^+_j(\l)\p _j^-(x, \l) -
m^-_j(\l)\p ^+_j(x, \l)\/ w_j(\l)}.
\]
 The function $e^{\mp ik_j(\l)x}\p ^{\pm}_{j}(x,\l)$ is
$\t_j-$periodic in $x$ for any $\l$, where $\t_1=\t, \t_2=1.$  Introduce the domain
$\L_j(\d)=\{\l\in \L_j, |\sqrt{\l}-\pi n/\t_j|\ge \d, n\ge 1\}, \d >0$.
Below we need the asymptotics
\[
\p ^{\pm}_{j}(x, \l)=e^{\pm i\sqrt{\l}x}(1+O(\l^{-1/2})),
\ \ \ \ m^{\pm}_{j}(\l)=\pm i\sqrt{\l} (1+O(\l^{-1/2})),
\ \ \ |\l | \to \iy , \
\]
$\l \in \L_j(\d)$, uniformly on $[0, \t_j]$ (see [T]).
We introduce $\P ^{\pm}(\cdot,\l)\in L^2({\R_{\pm}}), \l\in
\L_0^{(1)}=\C \sm \s_c(T_0)$,
which are solutions of Eq. (4.1) and having the form
$$
\P^{+}(x,\l)=\p_2^+(x,\l), \ \  x >0; \ \ \ \ \P^{-}(x,\l)=\p_1^-(x,\l), x <0,
\ \ \l\in \L_0^{(1)}.
$$
We determine these functions on the real line. First we find the Wronskians
\[
w(\l)=\{\P_-,\P_+\}=m_2^+(\l)-m_1^-(\l), \ \ \l\in \L_0^{(1)},
\]
\[
w_j(\l)=\{\p_j^-,\p_j^+\}=m_j^+(\l)-m_j^-(\l),\ \ \l\in \L_0^{(1)},\ \  j=1,2.
\]
Hence  we get
\[
w-w_2=(m_2^+-m_1^-)-(m_2^+-m_2^-) = m_2^--m_1^-,
\]
\[
w-w_1=(m_2^+-m_1^-)-(m_1^+-m_1^-) = m_2^+-m_1^+ .
\]
The definition of the Weyl functions yields
\[
w(\l)=\lt({a_2(\l)  \/ \vp_2(1, \l)}-{a_1(\l)  \/\vp_1(1, \l)}\rt) +
i\lt({\sin k_2(\l) \/ \vp_2(1, \l)}+{\sin k_1(\l) \/ \vp_1(1, \l)}\rt),
\ \ \ \l\in \L_0^{(1)}.
\]
Relations (4.6), (4.3) yield:
if $\l\in \s (H)\sm (\pa \s (H)\cup \{\m_n(p_j), j=1,2,\ n\ge 1 \})$, then
we get $\pm {\rm Im}\ w(\l\pm i0)>0 $. Using (4.1-2) we have
$$
\P ^+(x,\l)=\p^+_2(x,\l)=\vt_2(x,\l)+m_2^+(\l)\vp_2(x,\l)=v(x,\l)+m_2^+(\l)u(x,\l),
\ \ x \ge 0, \ \ \l\in \L_0^{(1)},
$$
and (4.2) implies $\P^+(x,\l)=\vt_1(x,\l)+m_2^+(\l)\vp_1(x,\l),\ \ \ x<0, \l\in\L_0.$
Then identities (4.3-9) yield for $x<0,\ \l\in \L_0^{(1)}$:
$$
\P^+={m_1^+\p _1^--m_1^-\p ^+_1 \/w_1}+m_2^+{\p^+_1-\p_1^-\/ w_1}=
{m_1^+-m_2^+ \/ w_1}\p _1^-+ {m_2^+-m_1^- \/ w_1}\p _1^+ =
{w_1-w \/ w_1}\p ^-_1+{w \/ w_1}\p _1^+,\ \ \
$$
and therefore,
\[
\P^+(x,\l)=\cases{\p_2^+(x,\l)\ \  &if\ \ \  $ x>0, $\cr
{w_1(\l)-w(\l)\/w_1(\l)}\p^-_1(x,\l)+{w(\l)\/w_1(\l)}\p_1^+(x,\l)\ \
                           &if\ \  $ x<0$ \cr},\ \ \ \l\in \L_0^{(1)},
\]
 and similarly we have
\[
\P^-(x,\l)=\cases{{w_2(\l)-w(\l)\/w_2(\l)}\p^+_2(x,\l)+{w(\l)\/w_2(\l)}\p_2^-(x,\l)
                 \ \ &if\ \ $ x > 0,$ \cr
\p_1^-(x,\l)\ \ \ \ \ &if\ \ $ x < 0$ \cr}, \ \ \ \l\in \L_0^{(1)}.
\]
Substituting (4.5) into (4.11-12) we obtain the asymptotics
\[
\P^{\pm} (x, \l)=\exp (\pm i\sqrt{\l}x)(1+O(\l^{-1/2})),
\ \ {\rm as}\ \ \ |\l | \to \iy , \ \l \in \L_1(\d)\cap \L_2(\d),
\]
uniformly on $[-\t,1]$, for some $\d>0$. Since the functions $\P_{\pm}(x,\l),
x\in\R$, are real on the half-line $\l\in (-\iy,0)$ we have the identities
\[
\bar\P^{\pm}(x,\l)=\P^{\pm} (x,\bar\l), \ \ \ \ \l\in \L_0^{(1)}.
\]
The kernel $R(x,x',\l)$ of the operator $(T_0-\l)^{-1}, \l\in\L$, has the form
\[
R(x, x',\l)={1\/w(\l)}\left\{\begin{array}{cc}
\P_+(x,\l)\P_-(x',\l)\ \ {\rm if}\ \ x>x'\\
\P_-(x,\l)\P_+(x',\l)\ \ {\rm if}\ \ x<x'\end{array}\right. , \ \ \ \l\in \L_0^{(1)}.
\]
We prove the first main theorem concerning  $T_t$.

\no {\bf Proof of Theorem 2.1} i) We shall use standard tools from
the theory of self-adjoint extensions and the Weyl Theorem. For
the sake of the reader, we shall briefly mention all arguments. We
need the following result of Krein [Kr]. Let $A$ be a closed
symmetric operator in a separable Hilbert space ${\cal H}.$ Define
the deficiency index $n_{\pm}=dim [ ker(A^{*}\pm i)]$. Remark that
by the Neumann extension theory, every closed symmetric operator
with $n_+=n_-$ admits a self-adjoint extension. In accordance with
Krein an open interval $(\a^-,\a^+),$ is called a gap of $A$ if
$$
\|(A-{\a^-+\a^+ \/ 2})f\|\ge {\a^+-\a^- \/ 2}\|f\|,\ \ \
 {\rm for \ \ \ \ all}\ \  \  \ f\in D(A).
$$
Assume that: a) $n_-=n_+=n\ge 0$, b) the operator $A$ has a gap $(\a^-, \a^+),$
c) there exists a self-adjoint extension $A_e$ of $A.$
Then for any self-adjoint extension $A_0$ of $A$ we have:

\no $\s _{ess}(A_0)=\s _{ess}(A_e),$  and the spectrum of $A_0$ is
discrete inside the gap and consists of at most $n$ eigenvalues
counting multiplicities.

We need the simple fact (see Zheludev [Z]). 
Define the operator $H^{\pm}_{m}y=-y''+p_my$ in $L^2(\R _{\pm})$, 
with the Dirichlet
boundary condition $y(0)=0$   for $m=1, 2,$. Furthermore,
$\s (H_m^{\pm})=\s_{ac}(H_m^{\pm})\cup \s_{d}(H_m^{\pm}), \
\s_{ac}(H_m^{\pm})=\s (H_m),$
and the eigenvalues of $H^{\pm}_{m}$ coincide with some $\m_n(p_m), n\ge 1$. Define the symmetric operator $Af=T_0f, \ f\in C_0^{\iy}(\R\sm\{0\})$. Then $n_-(A)=n_+(A)=2$ and there exists a
special self-adjoint extension $A_e=H^-_{1} \os H^+_{2}$ of $A$ in the space $L^2(\R_-)\oplus L^2(\R_+)$. 
Note that the operator $(T_0-i)^{-1}-(A_e-i)^{-1}$ is compact with rank 2, since $n_-(A)=n_+(A)=2$ (see [Kr]).

Hence using the results of Krein and Zheludev, (and the Weyl
Theorem about the invariance of the essential spectrum) we obtain:
 a) $\s _{ess}(T_0)=\s (H_1)\cup \s (H_2),$
\ b) there exist  gaps $\g _{nm}(T_0)=\g _{n}(H_1)\cap \g _{m}(H_2),
$ for some $n, m\ge 0,$  in the essential spectrum of $T_0$,\
c) for each $n, m\ge 0,$ the spectrum of $T_0$ is discrete
inside  the gap $\g _{nm}(T_0)\neq \es$ and consists of at most two eigenvalues.

Moreover, using the properties of $w, \P^{\pm} $ we deduce that the
function $R(x, x',\l\pm i0), \l\in \R ,$ is real analytic in $\l\in \R$ away from the discrete set of the zeros of the functions $w(\l),  \vp_m(1,\l), m=1,2$, which implies $\s _{sc}(T_0)=\es $.

We prove the absence of eigenvalues in the continuous spectrum
by contradiction. Let $\l\in \s_c (T_0), f $ be an eigenvalue and a corresponding eigenfunction of $T_0$. Hence $\l\in \s (H_1)$
or $\l\in \s (H_2)$. Let  $\l\in \s (H_2)$, the proof for 
$\l\in \s (H_1)$ is similar.
Then $-f''+p_2f=\l f, x>0.$ But it is well known that in this case
the function $f\not\in L^2(\R _+)$ (see [T] or [CL]). Hence
we have a contradiction.

\no The results of ii) will be proved in Lemma 4.1-3. 

iii) The Wronskian for the dislocation has the form
$w(\l,t)=m^+(\l,t)-m^+(\l,0)$, which follows from (4.10).
The Wronskian for the half-solid has the form
$w(\l,t)=m^+(\l,t)-\sqrt{s-\l}$, which follows from (4.10).
Using these identities and the definitions of $m^+(\l,t)$ (see(3.10)) we obtain the statement iii).
$\BBox$

\no   {\bf Lemma 4.1.}
{\it Let $q_0$ be a biperiodic potential in the sense of (1.1).
Then

\no  i) $\s^{(2)}=\s (H_1)\cap\s (H)\neq \es$ and for each $x,
x'\in \R $  the functions $R(x, x',\cdot ), w, w_1, w_2,
\P^{\pm}(x,\cdot )$ have a meromorphic continuation from
$\L_0^{(1)}$ across the set $\s^{(2)}$ to the second sheet
$\L_0^{(2)}$, where the following identities are fulfilled:
\[
w_n(\l^{(2)})=-w_n(\l), \ \ \ \ \ m^{\pm}_n(\l^{(2)})=m^{\mp}_n(\l),
\ \ \ \ \ \p ^{\pm}_n(\cdot , \l^{(2)})=\p ^{\mp}_n(\cdot , \l),\ \ \  n=1,2 ,
\]
\[
w(\l^{(2)})=m^-_2(\l)-m^+_1(\l),
\]
\[
\P^+(x,\l^{(2)})=\cases{\p_2^-(x,\l)\ \ \ &if\ \ \ \ $x>0, $\cr
                          \P^+(x,\l)+w_2(\l)\vp_1(x,\l) \ \ \ &if\ \ \ \ $x<0. $\cr}
\]
\[
\P^-(x,\l^{(2)})=\cases{\P^-(x,\l)+w_1(\l)\vp_2(v,\l)\ \ \ &if\ \ \ \ $x>0 $\cr
                     \p _1^+(x, \l) \ \ \ &if\ \ \ \ $x<0. $\cr}
\]
\no ii) The number $\l_e$ is an eigenvalue of $T_0 $ iff  the number $\l_e$
considered as a point on the second sheet $\l_e^{(2)}\in\L_0^{(2)}$ is a resonance
of $\wt T_0$. All resonances of $T_0$ on the sheet $\L_0^{(2)}$ lie on the gaps
$\g_n(T_0),  n\ge 0$ and $\#^{(2)}(T_0,\g_n)\le 2$ for all $n\ge 0$.

\no Proof.}  i) In order to get an analytic continuation we have to study the function
$\sin k_j(\l), j=1, 2$. The function $\sin k_j(\l)$  is real on
$\s (H_j)$. Then $k_j$ has a meromorphic continuation across $\s
(H_j)$ by the formula $\bar k_j(\l)=k_j(\bar \l)$. Using (3.7) we
have $k_j(\l)=\pi n+iv(\l), \l\in \g _n\ss\L_0^{(1)}$ and
$k_j(\l)=\pi n-iv(\l), \l\in \g_n\ss \L_0^{(2)}$, where $v=\Im k$.
Then we get $\sin k_j(\l)=-\sin k_j(\l^{(2)}), \l\in \g _n ,$ and
using (4.3), (4.6-7) we obtain (4.16-17); (4.16-17) and (4.11-12)
imply (4.18-19). Relations (4.16-19) yield a meromorphic
continuation of $R(x,x',\l)$ into the second sheet $\L_0^{(2)}$.

 ii) Let $\wt w$ be the Wronskian
for $\wt T_0$. Identity (4.17) yields $ w(\l^{(2)})=-\wt w(\l)$.
Then for any eigenvalue $\l_e$ of $T_0$ the number $\l_e^{(2)}$
is a resonance of $\wt T_0$ and inversely.
Hence all resonances on the sheet $\L_0^{(2)}$ of $T_0$  are real and lie
on the gaps $\g _n(T_0), n\ge 0,$ of this sheet.
Moreover, (1.3) yields $\#^{(2)}(T_0, \g _n)\le 2.$  $\BBox$  

We formulate the results about an analytic continuation
into the sheet $\L_0^{(3)}.$

\no    {\bf Lemma 4.2}
{\it Let $q_0$ be a biperiodic potential in the sense of (1.1)
and let $\s ^{(3)}=\s (H_1)\sm\s (H_2)\neq \es$. Then for each $x, x'\in \R $
the functions $R(x, x',\cdot ), w, w_1, w_2, \P ^{\pm}(x,\cdot )$ have meromorphic
continuations from  $\L_0^{(1)}$  across the set $\s ^{(3)}$
to the sheet $\L_0^{(3)}$, where the following identities are fulfilled:
\[
m^{\pm}_1(\l^{(3)})=m^{\mp}_1(\l),\ \ \ \p^{\pm}_1(\cdot,\l^{(3)})=
\p^{\mp}_1(\cdot,\l),\ \ \ w(\l^{(3)})=m^+_2(\l)-m^+_1(\l),
\]
\[
\P^-(x,\l^{(3)})=\cases {\P^-(x,\l)-w_1(\l)\vp_2^+(x,\l),\ \ \ &if\ \ \ \ $x>0,$\cr
                         \p_1^+(x,\l)\ \ \ &if\ \ \ \ $x<0, $\cr}
\]
$w_1(\l^{(3)})=-w_1(\l),$ and the functions $\P^+(x,\cdot ), w_2, m_2^{\pm}, \p_2^{\pm}$  don't change.}

We formulate the results about an analytical continuation into the sheet $\L^{(4)}.$

\no    {\bf Lemma 4.3} {\it 
Let $q_0$ be a biperiodic potential in the sense of (1.1)
 and let $\s^{(4)}=\s (H_2)\sm\s (H_1)\neq \es$. Then  for
each $x, x'\in \R $ the functions $R(x, x',\cdot ), w, w_1, w_2,
\P^{\pm}(x,\cdot )$ have meromorphic continuation from
$\L_0^{(1)}$ across the set $\s^{(4)}$ to  the sheet
$\L_0^{(4)}$ and the following identities are fulfilled:
\[
m^{\pm}_2(\l^{(3)})=m^{\mp}_2(\l),\ \ \p^{\pm}_2(x,\l^{(3)})=\p^{\mp}_2(x,\l),
\ \ w(\l^{(4)})=m^-_2(\l)-m^-_1(\l),
\]
\[
\P^+(x,\l^{(3)})=\cases {\p_2^-(x,\l)\ \ \ &if\ \ \ \ $x>0, $\cr
          \P^+(x,\l)-w_2(\l)\vp_1(x,\l)\ \ \ &if\ \ \ \ $x<0 $\cr}
\]
$w_2(\l^{(4)})=-w_2(\l), $ and the functions $\P^-(x,\cdot ), w_1, m_1^{\pm},\p_1^{\pm}$  don't change.}

The proof of Lemmas 4.2-3 repeats the proof of Lemma 4.1.

\section{Eigenvalues and resonances }
\setcounter{equation}{0}

We consider the eigenvalues of the operator $T_0$ with even
potentials $p_1, p.$ Recall that for an even potential the function 
$a(\cdot)\ev 0$ and the Weyl function $m^{\pm}=\pm i\sin k /\vp(1,\cdot)$.

\no    {\bf Proof of Theorem 2.2}
i) If $\wt\m_1<\wt\n_1$ and $\wt\m_2<\wt\n_2$, then the gaps have the form
$\wt\g_{1}=(\wt\m_1,\wt\n_1), \wt\g_{2}=(\wt\m_2,\wt\nu_2)$. Using (3.6-7), (3.24) we have $m_2^+<0$ and $m_1^->0$ in the gap $\wt\g (T_0)$. Then we deduce that
the Wronskian $w<0$ in this gap and $\#(T_0, \wt\g (T_0))=0.$

 Assume $\wt\mu _1>\wt\nu_1$ and $\wt\mu _2<\wt\nu_2$ and let $\wt\g (T_0)=(\a^-, \a^+)$.
The function $w$ is real analytic on $(\a^-, \a^+)$ and
relations (3.6-7), (3.24) imply $w(\a^-)<0, w(\a^+)>0$. Then the estimate  $\#(T_0, \wt\g (T_0))\le 2$ gives $\#(T_0, \wt\g(T_0))=1.$

\no ii) Using (3.6-7), (3.25)
we deduce that $m_1^->0$ on the interval $\wt\g (T_0).$ Let $\wt\mu _2<\wt\nu_2$.
Then (3.6-7), (3.25) yield $m_2^+<0$ on $(\wt\mu _2, \wt\nu_2)$,
 and the function $w<0$ on the interval $\wt\g (T_0).$
Let $\wt\n _2<\wt\m_2\le 0$.  Then $w(\wt\n_2)<0$  and $w(\wt\m_2)>0$. Then since
$\#(T_0, \wt\g (T_0))\le 2$, we get $\#(T_0, \wt\g (T_0))=1.$

\no iii) We consider the ground state in the gap $\wt\g _0(T_0)$.
Using (3.6-7), (3.24)  we obtain $m_2^+<0$ and $m_1^->0$
in the gap $\wt\g _0(T_0)$. Then the function $w<0$ on the interval
$\wt\g _0(T_0)$ and $\#(T_0, \wt\g _0(T_0))=0.$ $\BBox $

Next we obtain biperiodic potentials with a prescribed number of eigenvalues.

\no    {\bf Proof of Theorem 2.3.} Let $q_0\in Q$ and let $p\in
L_{even}^2(\T)$. Then due to Theorem 2.2 we have the sequence
$\{d_n\}_1^\iy$, $d_n=\# (T_0, \g_n(T_0))\in\{0,1\}, n\ge 1$. Note
that Theorem 2.1 yields $d_n=0$ if $|\g_n(T_0)|=0$. Hence we have
a mapping $q_0\to (p,d)\in P$.

\no Let $(p,d)\in P$. For fixed $p\in L_{even}^2(\T)$
the gaps have the forms: $\g_n(H)=(\m_n(p),\n_n(p))$ or $\g_n(H)=(\n_n(p),\m_n(p))$ for all $n\ge 1$ (see [GT]).
In order to construct $p_1$ we need the following result from [GT]:
for any $p\in L_{even}^2(\T)$ with gaps
$\g_n(H)=(\a_n^-,\a_n^+)$ 
there exists  unique $p_1\in L_{even}^2(\T)$ such that $\g_n(H)=\g_n(H_1)$ and
$\m_n(p_1)\in\{\a_n^-,\a_n^+\}$ for all $n\ge 1$.

In order to determine $p_1$, firstly, we set $\g_n(H)=\g_n(H_1)$.
Secondly, using the number $d_n\in \{0,1\}$, according to Theorem 2.2, we choose the position of $\m_n(p_1)$ by the following:
if $d_n=1$, then $\m_n(p_1)=\n_n(p)$ and if $d_n=0$, then $\m_n(p_1)=\m_n(p)$.
Thus, we know the gaps $\g_n(H_1)$ and the position of $\m_n(p_1)$ for all $n\ge 1$
and by the result of [GT], Theorem 2.2, there exists a unique isospectral $p_1\in L_{even}^2(\T)$ such that
$\#(T_0, \g_n(T_0))=d_n$ for all $n\ge 1$.
Hence the mapping $q_0\to (p,d)$ from $Q$ into $P$ is 1-to-1 and onto.

For isospectral potentials $p,p_1$ the identity $\|p_1\|=\|p\|$
holds (see [M] or [K1]) and the estimates (2.5) are fulfilled (see [K1]).
$\BBox$

We consider the eigenvalues for the general case. In order to
prove Theorem 5.1 we need the following definitions. Suppose that
a gap of $T_t$ has the form $\wt\g (T_t)=(\a^-, \a^+)=\g_n(H_2)\cap \g_m(H_1)\neq \es$
 for some gaps $\g_k(H_j)=(\a^-_{k,j}, \a^+_{k,j}), j=1,2, k=n,m$.
By Lemma 3.4, the equation $m_2^+(\a^-, y)=m_1^-(\a^-), y\in [0, 1]$
has $N\ge n$ roots $y_1, \dots , y_N$.
Assume for simplicity $y_1=0$ and introduce the function
\[
\cL (t,\a^-)=
\cases {\sqrt{2|M_{n,2}^-|}L_2(t,\a^-),\ \ \ &if \ \ \ \ $\a_{m,1}^-< \a_{n,2}^-, $\cr
\sqrt{2|M_{n,2}^-|}L_2(t,\a^-)/r(\a^-,t),\ \ \ &if\ \ \ \ $\a_{m,1}^-=\a_{n,2}^-, $\cr
\sqrt{2|M_{m,1}^-|}[\dot \P_1(0,\a^-)^2-(p_2(t)-\a^-)\P_1(0,\a^-)^2],
\ \ \ &if \ \ \ \ $\a_{m,1}^-> \a_{n,2}^-, $\cr}
\]
where 
$L_j(t,\l)=\mp [\dot \P_j(t,\l)^2-(p_j(t)-\l)\P_j(t,\l)^2]$ for $\l=\a_{k,j}^{\pm},\ k\ge 0, \ t\in [0, 1]$.
Here $\P_j(t,\a_{k,j}^{\pm})$  is the corresponding real normalized
eigenfunction for 2-periodic problem $-y''+p_jy=\l y$,
and $\pm M_{k,j}^{\pm}>0$ are the corresponding effective masses for $H_j, j=1,2$ and
$$
r(\l, t)=1+\sqrt{|M_{n,2}^-/M_{m,1}^-|}\P_2(t,\l)^2 /\P_1(t, \l)^2 ,
 \ \ \ \ \l=\a_{m,1}^-= \a_{n,2}^-.
$$
\no   {\bf Theorem 5.1.}  {\it
Let $q_t$ be a biperiodic potential  given by (1.1), where $p_1, p_2\in L^2_{loc}(\R)$.
Suppose that a gap $\wt\g (T_t)=(\a^-,\a^+)=\g_n(H_2)\cap \g_m(H_1)\neq \es$
 for some gaps $\g_k(H_j)=(\a^-_{k,j}, \a^+_{k,j}), j=1,2, k=n,m$.
Assume that $\m_{m,1}(p_1)\in (\a^-,\a_{m,1}^+]$ and let $m_2^+(\a^-, 0)=m_1^-(\a^-)$.
Then there exists a unique function $z(\cdot )\in W^2_1(-\ve, \ve),
z^{-}(0)=0$ for some $\ve >0$, such that

\no i) If $z(t)>0$, then $\l (t)\ev\a^-+z(t)^2\in\wt\g(T_t)$ is an eigenvalue of $T_t$.

\no If $z(t)<0$, then $\l (t)\in\wt\g (T_t)\ss\L_0^{(m_0)}$ is a resonance of $T_t$.
Here $m_0=2$ if $\a_{m,1}^-=\a_{n,2}^-$; \ $m_0=3$ if $\a_{m,1}^->\a_{n,2}^-$ and $m_0=4$ if
$\a_{m,1}^-<\a_{n,2}^-$. Moreover, the following asymptotics is fulfilled:
\[
z(t)=\int _0^t\cL (t,\a^-)dt+O(t^{3/2}) \ \ \ as \ \ t\to 0.
\]
ii)If $p_1, p_2$ are continuous then $z(\cdot )\in C^1(-\ve, \ve)$ and
$z(t)=\int _0^t\cL (t,\a^-)dt+O(t^2),  t\to 0$.

\no iii) If $p_1, p_2$ are analytic functions, then $z(\cdot )$
is real analytic on $(-\ve,\ve)$.

\no {\bf Remark}.  Result for an eigenvalue in some neighborhood
of $\a_n^+, n\ge 0$, is similar.

\no   Proof}. i) There are 3 cases: $\a^-_{n,2}<\a^-_{m,1}, a^-_{n,2}=a^-_{m,1}, a^-_{n,2}>a^-_1$.
Recall $ b_2(\l)=\sqrt {\D_2^2(\l)-1},$  and
\[
\z_2(\l,t)={\dot\vp_2(1,\l,t)\/2\vp_2(1,\l,t)},\ \
\ \ m_2^+(\l,t)=\z_2(\l,t)-{b_2(\l)\/\vp_2(1,\l,t)},
\ \ \ \l\in \g_2, \ \  t\in [0,1).
\]
Define a new parameter $z$ by the formula $z=z(\l)=\sqrt {\l-\a ^-}$, where
$z(\l)>0, \l\in \L$. Firstly, let $\a_{n,2}^-=\a^- > \a_{m,1}^-$.
By (3.9), the function $b_{02}(z)\ev b_2(\a^-+z^2)$ is analytic
in some neighborhood of zero and we get $b_{02}(z)=z\b_2^-(1+O(z^2)),$ as $ z\to 0$, where $\b_2^-\ev \D_2(\a^-)\sqrt{2|M_{n,2}^-|}>0$.
The function $1/\vp_2(1, \a^-+z^2, t)$ is smooth in the parameter
$(z, t)$ in some neighborhood of zero.  Therefore, $w(\a^-+z^2, t)$
is a smooth function of the parameter $(z, t)$ in some neighborhood
 of zero and $\dot w=\dot m^+_2$ (see (3.19)).  Moreover,
$w_z(\a^-+z^2, 0)|_{z=0}=-\b_2^-/\vp_2(1, \a^-, 0)\neq 0. $
Then by the implicit function Theorem there exists
a function $z(\cdot )\in W^2_1(-\ve , \ve )$
for some $\ve >0$ such that $w(\a^-+z^2(t), t)=0$ and $z(0)=0$.
We prove (5.2). Let $\O=1 +O(|z|+|t|)$. Relations (3.19-20), (3.9) yield
\[
\dot w(\l,t)=\dot m_2^+(\l,t)=\dot\z_2(\l,t)+O(z)=\dot\z_2(\a^-,t)+O(z),
 \ \ \ \ w_z(\l,t)=-{\b_2^-\O\/\vp_2(1, \a^-, t)},
\]
as $\ \ t\to 0$. Using (5.4), (3.22), we obtain
\[
\dot z(t)=-{\dot w(\l (t), t) \/ w_z(\l (t), t)}=
{\vp (1,\a^-,t)\dot \z_2(\a^-,t)+O(z)\/\b_2^-\O}=(\cL (t,\a^-)+O(z))\O.
\]
Let $z_m=\max |z(y)|, 0\le y\le t$. Then integration implies
\[
z(t)=\int _0^t \dot z(s)ds=
\int _0^t \cL (s,\a^-)ds+\sqrt{t}\|p_2\|O(|z_m|+|t|)+O(t^2+|tz_m|),
\]
hence $z(t)=O(\sqrt{t})$  and substituting this into (5.6) we get (5.2).

\no Secondly, let $\a_{n,2}^-=\a^- = \a_{m,1}^-$.
By (3.9), the function $b_{0m}(z)\ev b_m(\a^-+z^2), m=1, 2,$ is analytic in some neighborhood of zero and we get $b_{0m}(z)=z\b_m^-(1+O(z^2)), z\to 0$, where
$\b_m^-\ev\D_2(\a_-)\sqrt{|2M_{m,1}^-|}\neq 0$. The function $1/\vp_2(1, \a^-+z^2,t)$
is smooth in the parameter $(z, t)$ in some neighborhood of zero. Therefore,
$w(\a^-+z^2,t)$ is a smooth function of the parameter $(z,t)$ in some
neighborhood of zero.  Moreover,
\[
{\pa w(\a^-+z^2,0)\/\pa z}|_{z=0}=-{\b_2^-\/\vp_2(1,\a^-,0)}-
{\b_1^-\/\vp_1(1,\a^-)}=-{\b_2^-r(\a^-,0)\/\vp_2(1,\a^-,0)}\neq 0,
\]
since $\b_1^-\vp_2/(\b_2^-\vp_1)=(\b_2^-\P_2)^2/(\b_1^-\P_1)^2)$.
Then by the implicit function Theorem there exists a function
$z(\cdot )\in W^2_1(-\ve , \ve )$ for some $\ve >0$, such that
$w(\a^-+z^2(t), t)=0$ and $z(0)=0$. We find the asymptotics of $z(t)$ as  $t\to 0$. Using (3.19-20) we deduce that
\[
\dot w(\l, t)=\dot m_2^+(\l, t)=\dot \z_2(\l , t)+O(z)=
\dot \z_2(\a^-, t)+O(z), \ \ t\to 0,
\]
and by (5.7),
\[
w_z(\l_n(t), t)= \lt[{-\b_2^- \/ \vp_2(1, \a^-, 0)}+
{-\b_1^- \/ \vp_1(1, \a^-)}\rt]\O={-\b_2^-r(0,\a^-)\/\vp_2(1,\a^-,t)}\O .
\]
Then (5.8-9), (3.22) yield
\[
\dot z(t)=-{\dot w(\l (t),t)\/w_z(\l (t),t)}=
{\vp_2(1,\a^-,t)\dot\z_2(\l_n(t),t)+O(z)\/r(\a^-,0)\b_2^-\O}=(\cL (t)+O(z))\O.
\]
Then as above we get (5.2).

\no Next, we consider the last case $\a_{n,2}^-<\a^- = \a_{m,1}^-$. By (3.9), the function $b_{01}(z)\ev b_1(\a^-+z^2)$ is analytic in some neighborhood of zero and we get $b_{01}(z)=z\b_1^-(1+O(z^2)),$  as $ z\to 0.$ The function $1/\vp_2(\a^-+z^2, t)$ is smooth in the parameter $(z, t)$ in some neighborhood of zero.  Therefore, $w(\a^-+z^2, t)$ is a smooth function of the parameter $(z, t)$ in some neighborhood of zero. Moreover, $w_z(\a^-+z^2, 0)|_{z=0}=-\b_1^-/\vp_1(1,\a^-, 0)\neq 0.$  Then by the implicit
function Theorem there exists a function $z(\cdot )\in W^2_1(-\ve , \ve )$ for some $\ve >0$ such that $w(\a^-+z^2(t),t)=0$ and $z(0)=0$. We find the asymptotics of $z(t)$ as $t\to 0$. By (3.19-20),
\[
\dot w(\l,t)=\dot m_2^+(\l , t)=(p_2(t)-\l)-m_2^+(\l , t)^2=
(p_2(t)-\l)-\z_1^+(\a^-)^2+O(z), \ \ t\to 0,
\]
and by (3.9),
\[
w_z(\l , t)=-\b_1^- \vp_1(1,\a^-, 0)^{-1}\O.
\]
Hence the asymptotics (5.11-12) and the identity (3.22) imply
\[
\dot z(t)=-{\dot w(\l_n(t), t) \/  w_z(\l_n(t), t)} =
{\vp _1(1, \a^-)[(p_2(t)-\l)-\z_1(\a^-)^2+O(z)]\/\b_1^-\O}=(\cL (t,\a^-)+O(z))\O.
\]
Then, as above we get (5.2).

\no ii) If $p_, p_2$ is continuous, then using (5.5), (5.10), (5.13) we obtain
the needed estimates.

\no iii) If $p_1, p_2$ are analytic functions, then by the implicit function Theorem, $z(\cdot )$ is real analytic on $(-\ve,\ve )$. If $z(t)>0$ there exists an eigenvalue $\l(t)$ with the asymptotics (5.2). But if $z(t)<0$, then there is no an eigenvalue of $T_t$
 in the interval $[\a^-,\a^--\ve ]$ and we have a resonance.
$\BBox$

We study an eigenvalue $\l_0$ of $T_0$ in the 0-th gap $\g _0(T_0)$.

\no    {\bf Proposition 5.2}
{\it Let $q_0$ be a biperiodic potential  given by (1.1). Let $\g_0(H_j)=(-\iy,\a_{0,j}^+)$
be the infinite gaps of $H_j, j=1,2$ and $\min \{\a_{0,1}^+, \a_{0,2}^+\}=0$. Assume that
$\n_{0,j}\in\g_0(H_j)$ is the corresponding Neumann eigenvalue. Then $\# (T_0,(-\iy,\n^0])=0$,
where $\n^0=\min\{\nu_{0,1},\n_{0,2}\}$. Let in addition $w(0)>0,$ then $\#(T_0,\g_0(T_0))=1$.

\no Proof.} If $\l<\n^0 $, then by Lemma 3.5, $w(\l)=m_2^+(\l)-m_1^-(\l)<0$,
which yields the identity $\# (T,(-\iy ,\n^0])=0$. Let $w(0)>0.$  Then by (4.5),
$w(\l)\to -\iy $ as $\l\to -\iy $, and there exist
a zero of $w(\l), \l<0.$ But 
there is no other eigenvalue in the gap
since if we have 2, then we get 3 eigenvalues in the gap. This is impossible.
 $\BBox$

\section{Dislocations and half-solid}
\setcounter{equation}{0}

\no  {\bf 1. Dislocations.} We consider the operator
$T_t^{di}=-{d^2\/dx^2}+p_{(t)}$ acting in $ L^2(\R )$, where the dislocation potential $p_{(t)}=\c_-p+\c_+p(\cdot+t), t\in \R$ and $p\in L^1(\T)$. Recall that due to (2.6-7) we have only 
two sheets $\L_0^{(1)}, \L_0^{(2)}$ of the energy surface. Identities (3.10), (4.10) yield
\[
w(\l,t)=m^+(\l,t)-m^-(\l,0)={a(\l,t)-b(\l)\/\vp(1,\l,t)}-{a(\l,0)-b(\l)\/\vp(1,\l,0)}, \ \ \ \ \ \l\in \L_E(T_0).
\]
\no   {\bf Proof of Theorem 2.4.}
We consider the case $\a_n^-$ and let $n$ be even. The proof for $\a_n^+$ or for odd $n$ is similar. Let  $\m_n(t)=\m_n(p,t)$.
 Firstly, we consider the simple case $\a_n^-\neq \mu _n(0)$.
Define a new parameter $z$ by the formula $z=z(\l)=\sqrt {\l-\a ^-}$, where $z(\l)>0, \l\in \L_0$. Recall that the function $b_0(z)\ev b(\a_n^-+z^2)$
is analytic in some neighborhood of zero and has asymptotics (3.9).
The function $1/\vp(1, \a_n^-+z^2 , t)$ is smooth in the parameter
$(z, t)$ in some neighborhood of zero.  Therefore, $w(\a_n^-+z^2, t)$ is a smooth function of the parameter $(z, t)$ in some neighborhood of zero. Moreover, $w_z(\a_n^-+z^2, 0)|_{z=0}=-2\b_n/ \vp(1, \a_n^-, 0)\neq 0.$
Then by the implicit function Theorem there exists a function
$z(\cdot )\in W^2_1(-\ve , \ve )$ for some $\ve >0$ such that
$w(\a_n^-+z^2(t),t)=0, t\in (-\ve, \ve)$ and $z(0)=0$. We prove (2.8). By (3.19-20),
\[
\dot w(\l_n(t),t)=\dot m^+(\l_n(t),t)=\dot\z(\l_n(t),t)+O(z)=\dot\z(\a_n^-,t)+O(z),
\]
\[
w_z(\l_n(t),t)=-2\b_n/\vp(1,\a_n^-,t)\O, \ \ \O=1+O(|z|+|t|),
\]
as $t\to 0$. Recall $L(t,\l)=\dot \P(t,\l)^2-(p(t)-\l)\P(t,\l)^2$.
Using (6.2-3), (3.22) we obtain
\[
\dot z(t)=-{\dot w(\l_n(t),t)\/w_z(\l_n(t),t)}=
{\vp (1,\a_n^-,t)\dot\z(\l_n(t),t)+O(z)\/2\b_n\O}={\b_n\/2}L(t,\a_n^-)\O+O(z).
\]
Let $z_m=\max |z(y)|, 0\le y\le |t|$, then integration yields
\[
z(t)=\int _0^t\dot z(s)ds=
(\b_n/2)\int _0^t L(s, \a_n^-)ds+\|p\|\sqrt{t}O(|z_m|+|t|)+O(t^2+|tz_m|).
\]
Hence $z(t)=O(|t|^{1/2})$  and then $z_m= O(|t|^{1/2})$ and (6.6) yields (2.8).

 Secondly we consider the more complicated case
$\m_n(0)=\a_n^-$. Recall $\m_n(t)=\m_n(p,t)$. Let $a^0=a(\l,0), \vt_x^0=\vt_x(1,\l,0)$.
If we multiply (6.1) by $(a+b)$ and $(a^0-b)$ and use (3.11) we obtain
$-\vt_x (a^0-b)+\vt_x^0 (a+b)=0$. Therefore, we have the following
equation for the eigenvalues:
\[
\F(\l, t)\ev b(\vt_x+\vt_x^0)+(a\vt_x^0-a^0\vt_x)=0.
\]
Recall that $\vp (1,\l,t)=(\l-\m _n(t))\phi (\l ,t)$. Using (3.17) we have
\[
\dot \vp=-\dot \m_n\phi +(\l-\m_n)\dot \phi , \ \ \ \
\ddot \vp=-\ddot \m_n\phi -2\dot \m_n\dot \f+(\l-\m_n)\ddot \phi ,
\]
\[
\vp (1, \l, t)=O(|z|^2+t^2),\ \ \ \
\dot \vp (1, \l, t) =O(|z|^2+|t|),\ \ {\rm as}\ \ t\to 0,
\]
(6.8) yields $\F_z(\a_n^-,0)=2\b_n\vt_x(1,\a_n^-,0)\neq 0$ and the function
$\F(z,t)$ is analytic in $z$ and continuous in $t$. Then by the implicit function
Theorem, there exists a function $z(\cdot )\in W^2_1(-\ve,\ve )$ for some $\ve >0$
such that $\F(z(t), t)=0$ and  $z(0)=0.$ We will prove (2.9) and we have
\[
\F_z(z,t)=2\b_n\vt_x(1,\a_n^-,0)+O(|z|+|t|),\ \ \ \ {\rm as}\ \ \ t\to 0,
\]
and due to (3.12) we get
\[
\dot\F(z,t)=(\b-a^0)\dot\vt_x+\dot a\vt_x^0=
(\b-a^0)(\l-p)\dot\vp-\vt_x\vt_x^0-(\l-p)\vp\vt_x,
\]
\[
\dot\F(z,t)=-\vt_x(1,\a_n^-)^2+p(t)O(t^2+z^2)+O(|t|+|z|)
\]
Note that $\vt_x(\m_n(0),0)\neq 0$ and by (3.11),
$-(b^2)_{\l}= -\vt_x(1,\l,0)\vp_{\l}(1,\l,0), $ at $\l=\a_n^-$. Then (3.17) yields
\[
\vt_x(1,\a_n^-,0)=-2M_n^-/\vp_{\l}(1,\a_n^-,0)=\ddot \m_n(0)\f/2>0.
\]
Using (3.12), we have $\ddot \vp (1, \a_n^-, 0)=-2\b^2_n\dot \P(0, \a_n^-)^2$
and (6.9) yields $\ddot \vp (1, \a_n^-, 0)=-\ddot \m_n(0)\f $. Therefore, we get
${\vt_x^0/2\b_n}=(\b_n/2)\dot\P(0,\a_n^-)^2$. Using (6.8-12) we obtain
\[
z(t)=-\int_0^t{\dot\F(\l_n(t),t)\/\F_z(\l_n(t),t)}dt=
-t{\vt_x^0\/2\b_n}+\int_0^t[|p(t)|O(t^2+z^2)+O(|t|+|z|)]dt
\]
and since $p\in L^{2}(0, 1)$, (6.13) implies
\[
z(t)=t(\b_n/2)\dot\P(0,\a_n^-)^2+\sqrt{|t|}(\|p\|)O(|z_m|^2+t^2))+O(|z_mt|+t^2).
\]
Therefore, $z(t)=t(\b_n/2)\dot\P(0,\a_n^-)^2+\sqrt{|t|}O(|z_m|),$
and $z_m(t)=O(|t|^{1/2})$. Substituting these into (6.14) we have (2.9).
Note that $\l_n(t)=\a^-+z(t)^2$ while $z(\cdot )$ has the asymptotics (2.9)
and the function $a+b$ has a zero $\m_n(t)$ in the neighborhood of zero. It is important that $\l_n(t)$ and $\m_n(t)$  have different asymptotics, see (2.9) and (3.17).

Assume that $p$ is an analytic function.  Hence by the implicit
function theorem, there exists a unique analytic function $z_n^-(t)$
of the equation $w(t, z)=0$  for all $t\in (-\ve , \ve )$  for some
small $\ve >0.$ Then  $z(t)=rt^m(1+O(t)),$ for some $r\neq 0, m\ge 1.
\ \ \BBox$

We prove the existence of two eigenvalues in the gap for some potentials.

\no   {\bf Proof of Theorem 2.5.}
In Lemma 3.4 we have proved that for any finite sequences
$\{d_n\}^N, \{s_n\}^N, d_n\in\{0,1\}, s_n>0, N\ge 1$ there exists an even potential
$p\in C^2(\T)$  such that $\g_n(p)=s_n,$ and $L (0, \a_n^{\pm})>0,$
 if $d_n=1$ and  $L(0, \a_n^{\pm})<0,$ if $d_n=0$  for $n=1,..., N.$
We take the dislocation potential in the form  $p_{(t)}=\c_-p+\c_+p(\cdot+t)$.
Then by Theorem 2.4, there exists $\ve >0$ such that
$|\g _n(T_t^{di})|=s_n,$ and $\#(T_t^{di}, \g _n(T_t^{di}))=2d_n$ for all
$n=1, 2, .., N,$ and each $t\in (0, \ve )$.$\BBox$

\no  {\bf 2. Half-solid.}
We consider the half-solid operator $T_t^s=-{d^2\/dx^2}+q_t^s(x)$ acting in $L^2(\R )$, where the potential $q_t^s(x)=s\c_-+\c_+p(\cdot+t)$ and $s, t\in \R$,
$p\in L^1(\T)$ is real. Assume $\a_0^+(H)=0$. We have the following simple results concerning $\s (T_t^s).$  If $s\le \a_1^-(H)$, then $\s_{ac}(T_t^s)= (s_-,\iy )$, where $s_-=\min (0,s)$. If $s>\a_1^-$ then there is a gap in the spectrum of $T_t^s.$ Our goal is to study the eigenvalues in the gaps $\g_n(T_t^s), n\ge 0,$
and to find how these eigenvalues depend on $t, s.$ It is clear that they depend on $t$ periodically. Using (4.6) we rewrite the Wronskian $w(z,t)$ in the form
\[
w(\l,t)=m^+(\l,t)-\sqrt{s-\l}={a(\l,t)-b(\l)\/\vp(1,\l,t)}-\sqrt{s-\l}, \ \ \ \l<s,
\]
where $\sqrt{s-\l}>0$ if  $\l<s, \l\in \L_0^{(1)}$.

\no {\bf Proof of Theorem 2.6}. i) Using Theorem 2.2 with $p_1=s, p_2=p$,
we obtain  (2.12) and the identity $\#(T_0^s, \g _0(T_0^s))=0.$  By Lemma 4.3,
the Wronskian $w$ on the second sheet $\L^{(4)}$ has the form
$w(\l^{(4)})=b(\l)\vp(1,\l)^{-1} -\sqrt{s-\l}$. Repeating the proof of Theorem 2.2 we get $\#^{(4)}(T_0^s, \g _n(T_0^s))=1-\#(T_0^s, \g _n(T_0^s))$. Moreover, (2.4) yields (2.12).

\no ii) Let $q_0^s\in Q_N, N\ge 0$ and $p\in L^2(\T)$ be even. Then
due to Theorem 2.2 we have the sequence $\{d_n\}_1^N$, $d_n=\#
(T_0^s,\g_n(T_0^s)), n\ge 1$, and $\ve \in (0,1)$ and $r\in
\ell^2$. Note that Theorem 2.1 yields $d_n=0$ if $|\g_n(T_0)|=0$.
Hence we have a mapping $q_0^s\to (r,d,\ve )\in P$.

Let $(r,d,\ve )\in P_N, N\ge 0$. Define the signed gap length
$l_n(p)=(\m_n(p),\n_n(p)), n\ge 1$ (see [GT]). For the sequence $r\in \ell^2$, there exists an even potential $p\in L^2(\T)$ such that $|l_n(p)|=r_n\ge 0, n=1,.., N,  \ l_n(p)=r_n, n\ge N+1.$
In order to get uniqueness of $p$ we need to determine $\sign l_n(p), n=1,.., N$, (see [GT]).

\no Fix a number $n=1,.., N$ and the gap $\g_n(H)=(\a_n^-,\a_n^+)\neq \es$.
We have 2 cases. First, if $d_n=1$, then using (2.12) we obtain
$l_n(p)>0$. Second, if $d_n=0$, then using (2.12) we obtain
$l_n(p)<0$. Hence we have $\{l_n\}\in \ell^2$ and there
exists a unique $p\in L^2(\T)$. Using $\ve \in (0,1)$ and the
definition of $\ve$ we obtain $s$, which yields a unique $q_0^s\in Q_N$.
\ \ \ $\BBox$

By Lemma 3.4, for $n\ge 1$, Eq. (2.13) has $N\ge n$ roots $y_1, y_2,\dots, y_N \in [0,1)$. The case $n\ge 1$ arises for the gap $\g _n(T_t(s)), n\ge 1,$  and $n=0$ corresponds the basic gap $\g _0(T_t^s)$.

\no {\bf Proof of Theorem 2.7} i)
We consider the case $\a_n^-< s$. The proof for  $\a_n^+< s$
is similar. The proof follows from i) of Theorem 5.1.
In this case the Weyl function has the form (6.16)
and $\cL (y,\l)=-[\P'(y,\l)^2+(\l-p(y))\P (y,\l)^2], \l=\a_n^-$.
The identity (3.5) yields $\z(y, \a_n^-)=\P '(y, \a_n^-)/\P (y, \a_n^-)$.
Together with (6.16) this implies $w(\a_n^-, y)=0.$ All conditions of
Theorem 5.1, i) are fulfilled and we get $z_n^{\pm}(t)$ with the needed properties. By (2.13), $\P_y(y,\l )=\sqrt{s-\l}\P (y,\l )$, then  $2\cL (y,\a_n^-)=(p(y)-s)\P (y,\a_n^-)^2$ and (5.2) yields (2.14).

\no ii) In this case we use i) of Theorem 5.1, when $s$ is the end of the gap $(\a_n^-,s)$ of the operator $T_y(s)$. We have
 $\P_1\ev 1, \P_1'\ev 0, M_1^+=1/2$, then  $2\cL (y,s)=(p(y)-s)$
and Theorem 5.1 (see the remark after this Theorem) implies the needed results including (2.15). $\BBox$

Now we consider eigenvalues (the ground state) in the gap $\g _0(T_t^s)$. The existence of eigenvalue is connected with some property of the function
$m^+(0,t),t\in[0,1]$, when $b>0, \vp>0$ in the gap $\g _0(T_t^s)$. We have the following result.

\no    {\bf Theorem 6.1}
{\it Let $T_0^s=-{d^2\/dx^2} +s\chi_-(x)+\chi_+(x)p(x)$,
where $p\in L^1(\T)$ is real and let $\a_0^+=0$.
If $m^+(0)<0 $, then  $\#( T_0^s, \g _0(T_0^s))=0 $ and if $ m^+(0)>0$, then
\[
\#(T_0^s,\g_0(T_0^s))=\cases{0, \ \ &if \ \  $s\le \n_0\ {\rm or}\ s\ge m^+(0)^2$\cr
                             1, \ \ &if \ \ $\nu_0<s< m^+(0)^2$\cr}.
\]
\no Proof.} Let $m^+(0)<0 $, then Lemma 3.5 yields $w(\l)<0$ for $\l<\min\{s,0\}$.
Hence eigenvalues in $\g_0(T_0^s)$ are absent.

\no Let $m^+(0)>0 $. Assume $\l<s\le \n_0$.  Again by Lemma 3.5, $m^+(\l)<0$ for
$\l< \n_0$. Then $w(\l)<0$  for $\l<s\le \n_0$, and eigenvalues  in $\g_0(T_0^s)$
are absent. Now, assume $m^+(0)^2\le s.$  Let the function $w$ have a simple zero
 $\l_1\in (\nu_0, 0)$. Then for some $s_0>s$ the function $m^+(\l)-\sqrt{s_0-\l}$
has a multiple root  $\l_1\in (\nu_0, 0)$. But this is impossible. Then
the first line in (6.17) is true.
 Consider the case $\nu_0<s<m^+(0)^2.$ Then by Lemma 3.5, $w(\nu_0)<0$ and
$w(\l_1)>0, \l_1=\min \{s,0\}$. Hence Proposition 5.2 yields the second line in (6.17).
$\BBox$

\no {\bf References}

\no [A] Anoshchenko, O.: The inverse scattering problem for the
Schr\"odinger equation with a potential that has periodic
   asymptotics.  J. Soviet Math. 48 (1990), no. 6, 662--668

\no [A1]  Anoshchenko, O.: Eigenfunction expansion of the
Schr\"odinger equation with a potential that has periodic
   asymptotics. J. Soviet Math. 49 (1990),   no. 6, 1237-1241

\no [BS] Bikbaev R., Sharipov R.: Asymptotics as $t\to \iy $ of
the Cauchy problem for the Korteveg-de Vries equation in the class
of potentials with finite-gap behavior at $x\to \pm \iy $,
Theoret. and Math. Phys. 78, 1989.

\no [CL]  Coddington, E., Levinson, N.:
Theory of ordinary differential equations. McGraw-Hill Book Company, Inc.,   New York-Toronto-London, 1955.

\no [DS] Davies E., Simon B. Scattering theory for systems with different spatial  asymptotics on the left and right, Commun. Math. Phys. 63, 277-301, 1978

\no [GT] Garnett J., Trubowitz E.: Gaps and bands of one dimensional
 periodic Schr\"odinger operator. Comment. Math. Helv. 59, 258-312 (1984).

\no [KK1] Kargaev P., Korotyaev E.: The Inverse Problem for
the Hill Operator, the Direct Approach, Invent. Math., 129, no. 3, 567-593 (1997)

\no [KK2] Kargaev P., Korotyaev E.: Effective masses
and conformal mappings. Commum. Math. Phys. 169, 597-625 (1995).

\no [K1] Korotyaev E.:  Estimates for the Hill operator, I.
J. Diff. Eq. 162, 2000, 1-26

\no [K2] Korotyaev E.: Inverse Problem and the trace formula for
the Hill Operator, II,    Math. Z., 231, 345-368 (1999)

\no [K3] Korotyaev E.:  Lattice dislocations in 1-dimensional model.
Commun. Math. Phys. 213, 471-489, 2000

\no [KP] Korotyaev E., Pokrovski, A.: One dimensional  half - crystal (in preparation)

\no [Kr] Krein M.: Theory of selfadjoint extensions of semi-
bounded Hermitian operators and its applications. I. Mat. Sbornik
20, No. 1, 1-95, 1947.(Russian)

\no [L] Levitan B.M.: The Inverse  Sturm-Liouville Problems.
Moscow, Nauka, 1984.(Russian)

\no [M] Marchenko V.: Sturm-Liouville operator and applications.
Basel: Birkh\"auser 1986.

\no [Mos] Moser J. : An Example
of a Scr\"odinger operator with almost periodic potential and
nowhere dense spectrum. Comment. Math. Helv. 56, 198-224 (1981).

\no [Ta] Tamm I. Phys.Z. Sowjet.1, 1932,  733.

\no [T] Titchmarsh E.: Eigenfunction expansions associated with
second-order differential equations 2, Clerandon Press, Oxford,
1958.

\no [Tr] Trubowitz E.: The Inverse problem for Periodic
Potentials . Commun. on Pure and Applied Math. V. 30, 321-337,  1977.

\no  [PTr] P\"oschel, P.; Trubowitz E.: Inverse Spectral Theory.
Boston: Academic Press, 1987.

\no [Z] Zheludev V. On the spectrum of Schr\"odinger operator with
periodic potentials on the half-line. Trudy kafedry  mat. anal.
Kaliningarad Univ. 1969.

\end{document}